\documentclass[reqno, english,11pt]{smfart}

\usepackage{amsfonts, amsthm, amsmath, amssymb }
\usepackage{a4wide}

\RequirePackage{mathrsfs} \let\mathcal\mathscr

\numberwithin{equation}{section}

\newtheorem{theorem}{Theorem}[section]
\newtheorem{lemma}[theorem]{Lemma}

\theoremstyle{definition}
\newtheorem*{ack}{Acknowledgements}

\newcommand{\oh}{\mathfrak o}

\newcommand{\No}{N_{W}(F,P)}

\newcommand{\Bl}{B^{(l)}}
\newcommand{\Blp}{B^{(l)}_1}

\newcommand{\wtil}{\mathfrak{I}}

\newcommand{\bbb}{\mathfrak{b}}
\newcommand{\sigstar}{\sideset{}{^*}\sum_{\sigma\bmod{\mathfrak{b}}}}
\newcommand{\ord}{\mathrm{ord}}
\newcommand{\starsum}{\sideset{}{^*}\sum}
\newcommand{\sumtwo}{\sideset{}{^{(2)}}\sum}

\newcommand{\fbth}{\fb_3}

\renewcommand{\d}{\mathrm{d}}
\renewcommand{\phi}{\varphi}

\newcommand{\0}{\mathbf{0}}

\newcommand{\PP}{\mathbb{P}}
\renewcommand{\AA}{\mathbb{A}}

\newcommand{\ZZ}{\mathbb{Z}}

\newcommand{\QQ}{\mathbb{Q}}
\newcommand{\RR}{\mathbb{R}}
\newcommand{\CC}{\mathbb{C}}

\newcommand{\cM}{\mathcal{M}}
\newcommand{\cW}{\mathcal{W}}

\renewcommand{\leq}{\leqslant}

\renewcommand{\geq}{\geqslant}

\newcommand{\ma}{\mathbf}

\newcommand{\m}{\mathbf{m}}
\newcommand{\M}{\mathbf{M}}
\newcommand{\x}{\mathbf{x}}
\newcommand{\y}{\mathbf{y}}
\renewcommand{\c}{\mathbf{c}}
\newcommand{\f}{\mathbf{f}}
\renewcommand{\v}{\mathbf{v}}
\renewcommand{\u}{\mathbf{u}}
\newcommand{\z}{\mathbf{z}}
\renewcommand{\b}{\mathbf{b}}
\renewcommand{\a}{\mathbf{a}}
\renewcommand{\k}{\mathbf{k}}
\newcommand{\h}{\mathbf{h}}
\newcommand{\g}{\mathbf{g}}
\renewcommand{\j}{\mathbf{j}}
\renewcommand{\r}{\mathbf{r}}

\newcommand{\fo}{\mathfrak{o}}
\newcommand{\fa}{\mathfrak{a}}
\newcommand{\fb}{\mathfrak{b}}
\newcommand{\fc}{\mathfrak{c}}
\newcommand{\fd}{\mathfrak{d}}
\newcommand{\fg}{\mathfrak{g}}
\newcommand{\fh}{\mathfrak{h}}

\newcommand{\fp}{\mathfrak{p}}
\newcommand{\fq}{\mathfrak{q}}
\newcommand{\fz}{\mathfrak{z}}
\newcommand{\B}{\mathbf{B}}

\newcommand{\ve}{\varepsilon}
\newcommand{\e}{\ensuremath{\mathrm e}}

\newcommand{\bbe}{\boldsymbol{\beta}}
\newcommand{\brh}{\boldsymbol{\rho}}
\newcommand{\bxi}{\boldsymbol{\xi}}

\DeclareMathOperator{\supp}{supp}

\DeclareMathOperator{\Res}{Residue}
\DeclareMathOperator{\tr}{Tr}
\DeclareMathOperator{\nm}{Nm}
\DeclareMathOperator{\n}{N}
\newcommand{\fS}{\mathfrak{S}}
\newcommand{\cH}{\mathcal{H}}

\DeclareMathOperator{\Mod}{mod} 
\renewcommand{\bmod}[1]{\,(\Mod{#1})}

\begin{document}
\title[The circle method for number fields]{Cubic hypersurfaces and a version of the circle method for number fields}
\author{T.D. Browning}
\address{School of Mathematics\\
University of Bristol\\ Bristol\\ BS8 1TW
\\ United Kingdom}
\email{t.d.browning@bristol.ac.uk}

\author{P. Vishe}
\address{
Department of Mathematics\\
University of York\\
York\\ YO10 5DD\\ United Kingdom}
\email{pankaj.vishe@york.ac.uk}

\begin{abstract}
A version of the Hardy--Littlewood circle method is developed for number fields 
$K/\QQ$ and is used to show that non-singular projective cubic hypersurfaces over $K$ always have a $K$-rational point when they have dimension at least $8$.
\end{abstract}


\maketitle

\setcounter{tocdepth}{1}
\tableofcontents

\section{Introduction}

Much of analytic number theory has only been extensively developed for problems defined over the rational numbers $\QQ$.   While many generalisations to finite extensions of $\QQ$ are straightforward, there remain a number of  areas where substantial technical obstructions persist. One such lacuna may be found in the Hardy--Littlewood circle method,  which over $\QQ$ begins with a generating function 
$$
S(z)=\sum_{n\in \ZZ} a(n) e^{2\pi i z n}.
$$
Typically one is interested in the term $a(0)$,  detected via
$$
a(0)=\int_0^1 S(z) \d z,
$$
 the idea being to break $[0,1]$ into subintervals $[a/q-\delta,a/q+\delta]$ 
over which the integral is more easily estimated.
A key innovation, due to Kloosterman \cite{K}, involves 
decomposing $[0,1]$ using 
a Farey dissection
to keep track of the precise endpoints of the intervals. This allows one to introduce non-trivial averaging over the numerators of the approximating fractions $a/q$, an approach that is usually called the ``Kloosterman refinement''.  This method is not immediately 
available to us when passing to finite extensions of $\QQ$,
since 
aside from work of Cassels, Ledermann and Mahler \cite{cassels} particular to the imaginary quadratic fields $\QQ(i)$ and $\QQ(\rho)$, 
no  generalisation is known of the Farey dissection 
to the number field analogue of $[0,1]$.

The primary aim of this paper is to provide an alternative route to the Kloosterman refinement over an arbitrary number field,  circumventing the need for a Farey dissection.  We will illustrate the utility 
of this new approach by applying it to a long-standing problem in Diophantine geometry. 
Given a  
cubic hypersurface $X\subseteq \PP_K^{n-1}$ defined over a number field $K$, a ``folklore'' conjecture predicts that the set $X(K)$ of $K$-rational points on $X$ is non-empty as soon as $n\geq 10$. 
The following result establishes this conjecture for generic cubic hypersurfaces.

\begin{theorem}\label{t:10}
Let $K$ be a number field and let $X\subseteq \PP_K^{n-1}$ be a non-singular cubic hypersurface defined over $K$. If $n\geq 10$ then $X(K)\neq \emptyset$.
\end{theorem}

In fact, as conjectured by Colliot-Th\'el\`ene \cite[Appendix A]{PV}, we expect the Hasse principle to hold for non-singular cubic hypersurfaces $X\subseteq \PP_K^{n-1}$ with $n\geq 5$.  Work of Lewis \cite{lewis} ensures that  
$X(K_v) \neq \emptyset$ for every valuation $v$ of $K$ when $n\geq 10$. Hence  Theorem \ref{t:10} confirms the Hasse principle for non-singular cubic hypersurfaces  in $n\geq 10$ variables.  

The resolution of Theorem \ref{t:10} for the case $K=\QQ$ goes back to groundbreaking work of Heath-Brown \cite{hb-10}. 
Extending this approach to general  number fields $K$,
the best result in the literature is due to Skinner \cite{skinner}, who  requires $n\geq 13$ variables. The loss of precision is entirely due to the lack of a suitable Kloosterman methodology, 
a situation that we remedy in the present investigation. 
When no constraints are placed on the singular locus of $X$, work of Pleasants \cite{pleasants} shows that $n\geq 16$ variables are needed to ensure that $X(K)$ is non-empty. 
Finally, if the singular locus of $X$ contains a set of three conjugate points then Colliot-Th\'el\`ene and Salberger \cite{plus} have shown that the Hasse principle holds provided only that $n\geq 3$. 

It is now time to present the main technical tool  in this work.  
Let  $K$ be a number field of degree $d$ over $\QQ$, with ring of integers $\fo$.
The  ideal norm will be designated $\n \fa= \#\fo/\fa$ for any integral ideal
$\fa\subseteq \fo$.  In line with our description of the Hardy--Littlewood circle method, 
we would like to use  Fourier analysis to detect 
when elements of $\fo$ are zero. In fact we will be able to handle 
the indicator function
$$
\delta_K(\fa) =
\begin{cases}
1,  & \mbox{if $\fa = (0)$}, \\
0,  & \mbox{otherwise},
\end{cases}
$$
defined on integral ideals $\fa\subseteq \fo$.
When $K=\QQ$ an extremely useful formula for $\delta_\QQ$ 
was developed by Duke, Friedlander and
Iwaniec \cite{DFI}. This was later revisited by Heath-Brown \cite[Thm.~1]{H} in an effort to
relate it to classical versions of the circle method. In this paper we adapt the
latter approach to the setting of arbitrary number fields $K$,  as follows.

\begin{theorem}
\label{t:delta}
Let $Q\geq 1$ and let $\fa\subseteq \fo$ be an ideal.
Then there exists a positive constant $c_Q$ and an infinitely differentiable function $h(x,y): 
(0,\infty)\times \mathbb R\rightarrow \RR$ such that
$$
 \delta_K(\fa)=\frac{c_Q}{Q^{2d}}
 \sum_{(0)\neq \mathfrak{b}\subseteq \fo  }
~\sigstar \sigma(\fa)h\left(\frac{\n \fb}{Q^{d}}  , \frac{\n \fa}{Q^{2d}}\right),
$$
where
the notation 
$\sum^{*}_{\sigma \bmod{\fb}}$ means that the sum is taken over
primitive additive characters modulo $\fb$ extended to ideals, as 
described in \S \ref{s:characters}.
The constant $c_Q$ satisfies 
$$
c_Q=1+O_N(Q^{-N}),
$$ 
for any  $N>0$. Furthermore,  we have $h(x,y)\ll x^{-1}$ for all $y$
and $h(x,y)\neq 0$ only if $x\leq \max\{1,2|y|\} $.
\end{theorem}

This result will be established in \S
\ref{s:delta}.  
The number field $K$ is considered fixed once and for all. Thus 
all implied constants in our work are allowed to
depend implicitly on $K$. 
One obtains a formula for the indicator function on $\fo$ by restricting to principal ideals, in which case one writes $\delta_K((\alpha))=\delta_K(\alpha)$ for any $\alpha\in \fo$.

Given the broad impact that  \cite{DFI} and \cite{H} have had on number-theoretic problems over $\QQ$, one might 
view Theorem \ref{t:delta} as foundational in a systematic programme of work to 
extend our understanding to the setting of general number fields. 
In \S \ref{s:poisson-forms} we will indicate how Theorem~\ref{t:delta} can be used to count suitably constrained 
$\fo$-points on hypersurfaces. 
The outcome of this is recorded in Theorem \ref{prop:2}.
Once applied to cubic forms this will serve as the footing
for our proof of Theorem \ref{t:10}.
Furthermore, although we will not give details,  it can be applied analogously to study the density of $\fo$-points on 
hypersurfaces defined by 
quadratic polynomials $Q(X_1,\ldots,X_n)-m$, with $Q$ a non-singular quadratic form defined over $\fo$ and  $m\in \fo$
non-zero.
When $n\geq 5$ this is covered by work of Skinner \cite{skinner2}. 
Handling the case $n=4$, however,  makes essential use of Theorem \ref{t:delta}. Indeed, 
when $K=\QQ$, 
 it was precisely in this context that  Kloosterman's method originally arose.
Feeding this into the strategy of Eskin, Rudnick and Sarnak \cite{ERS}, one could use this result to give a new proof of Siegel's mass 
formula over general  number fields.

\begin{ack}
This work was initiated during the programme
``Group actions in number theory'' at the
{\em Centre Interfacultaire Bernoulli}
 in Lausanne, 
the  hospitality and financial
support of which is gratefully acknowledged. 
The second author would like to thank Akshay Venkatesh for introducing him to the problem and subsequent discussions, in addition to Brian Conrad and Pieter Moree for 
 helpful conversations.
The authors are very grateful to the anonymous referee for 
numerous pertinent comments.
While working on this paper the first author was 
supported by ERC grant 
 306457 and the second author was partly supported by EPFL,
  the G\"oran Gustafsson Foundation (KVA) at KTH and MPIM.
\end{ack}

\section{Technical preliminaries}
\label{s:technical}

Our work will require a good deal of notation. In this section we collect together the necessary
conventions, in addition to some preliminary technical tools, relevant to our number field 
 $K$ of  degree $d$ over $\QQ$.
Let $r_1$ (resp.\ $2r_2$) be the  number of distinct real (resp.\ complex) embeddings of $K$, with $d=r_1+2r_2$.
Given any $\alpha \in K$ we will denote the norm and trace by
$\n_{K/\QQ}(\alpha)$ and $\tr_{K/\QQ}(\alpha)$, respectively.  
Let $\rho_{1},\dots,\rho_{r_1}$
be the $r_1$ distinct real embeddings  and let 
$\rho_{r_1+1},\dots , \rho_{r_1+2r_2}$ be  a complete  set of $2r_2$ 
distinct complex embeddings, with 
$\rho_{r_1+i}$ conjugate to $\rho_{r_1+r_2+i}$ for $1\leq i\leq r_2$.  
Let $V$ denote the  $d$-dimensional commutative 
$\RR$-algebra
$$
K\otimes_\QQ \RR 
\cong \bigoplus_{l=1}^{r_1+r_2} K_l,
$$
where
$K_l$ is the completion of $K$ with respect to $\rho_l$, for $1\leq l\leq r_1+r_2$.
Thus $K_l=\RR$ (resp.\ $K_l=\CC$) for $1\leq l\leq r_1$ (resp.\  $r_1<l\leq r_2$).
Given  $v=(v^{(1)},  \dots, v^{(r_1+r_2)}) \in V$ we define 
\begin{equation}\label{eq:norm-trace}
\begin{split}
\nm(v) &= 
v^{(1)}\cdots v^{(r_1)}|v^{(r_1+1)}|^{2}\cdots |v^{(r_1+r_2)}|^{2}, \\
\tr(v)
&= v^{(1)}+\cdots+v^{(r_1)}+2\Re  (v^{(r_1+1)})+\cdots + 2\Re  (v^{(r_1+r_2)}).
\end{split}
\end{equation}
Furthermore, we define the character 
$
\e(\cdot)=e^{2\pi i \tr (\cdot)}
$
on  $V$. 
We will typically write $v^{(l)}\in K_l$ for the
projection of any $v\in V$ onto the $l$th component, for $1\leq l\leq r_1+r_2$. Thus 
any $v\in V$ can be written $v=\bigoplus_l v^{(l)}$.
Likewise, given a vector $\v\in V^n$, we will usually denote by $\v^{(l)}\in K_l^n$ the projection of the vector onto
the $l$th component.

There is a canonical embedding of $K$ into $V$ 
given by $\alpha
\mapsto  (\rho_{1}(\alpha), \dots ,  
\rho_{r_1+r_2}(\alpha))$, 
and we shall identify $K$ with its image in $V$.
Under this identification  any fractional ideal becomes a lattice in $V$. 
Let  $\{\omega_{1},\dots,\omega_{d}\}$ be a $\ZZ$-basis for $\fo$.  Then
 $\{\omega_{1},\dots,\omega_{d}\}$  forms an $\RR$-basis for $V$ and we may view 
$V$ as the set 
$\{x_{1}\omega_{1}+\cdots +x_{d}\omega_{d}: x_{i}\in \RR\}$.

We will need to introduce some norms on $V$ and $V^n$.
To begin with let  
$$
\langle v \rangle = \max_{1\leq l\leq r_1+r_2}|v^{(l)}|,
$$
for any $v\in V$.  We extend this to $V^n$ 
in the obvious way.
Next, let
\begin{equation}\label{eq:cl}
c_l=
\begin{cases}
		1,  & \mbox{if  $1\leq l\leq r_1,$} \\
		2, & \mbox{if  $r_1<l \leq r_1+r_2$.}
\end{cases}
\end{equation}
We will also need to
introduce a Euclidean norm $\|\cdot \|$ on $V$, given by 
$$
\|v\|= \sqrt{\sum_{1\leq l\leq r_1+r_2} c_l|v^{(l)}|^2
}.
$$
We extend this to $V^n$ by setting 
$$
\|\v\|= \sqrt{\sum_{1\leq i\leq n} \|v_i\|^2
}=
 \sqrt{\sum_{1\leq l\leq r_1+r_2} c_l|\v^{(l)}|^2
},
$$
if $\v=(v_1,\ldots,v_{n})\in V^n$.

We will make frequent use of  the dual form
with respect to the trace. For any fractional 
ideal $\mathfrak{a}$ in $K$ one defines the dual ideal
$$
\hat \fa = \{\alpha\in K: \mbox{$\tr_{K/\QQ}(\alpha x)\in  \ZZ$ for all $x\in \mathfrak{a}$}\}.
$$
In particular 
$\hat \fa= \fa^{-1}\fd^{-1}$,  where 
$$
\fd = 
\{\alpha\in K: \alpha\hat \fo\subseteq \fo\}
$$ 
denotes the different ideal of $K$ and is itself an integral ideal. One notes that $\hat \fo=\fd^{-1}$. Furthermore,
we have $\hat \fa\subseteq \hat \fb$ if and only if $\fb\subseteq \fa$.
An additional integral ideal featuring in our work is the denominator ideal
$$
\fa_\gamma=\{ \alpha\in \fo: \alpha \gamma\in \fo\},
$$
associated to   any $\gamma\in K$.

We let $D_K$ denote the modulus of the discriminant of $K$. 
Finally, we will reserve $U_K$ 
for denoting the set of units in $\fo$. For an integral ideal
$\fa$,  the notation $\sum_{\alpha\in \fa/U_K}$  means that
the  sum is over elements in $\fa$ modulo the action of $U_K$. 
Given an element $v\in V$ it will sometimes prove  advantageous  to use the action of $U_K$ to control the size of each component $v^{(l)}$ of $v$. 
This is the object of the following standard result. 

\begin{lemma}\label{lem:unit}
Let $v=\bigoplus_l v^{(l)}\in V$.  Then there exists $u\in U_K$ such that 
$$
|\nm(v)|^{1/d}\ll |( u v)^{(l)}|\ll |\nm(v)|^{1/d},
$$
for $1\leq l\leq r_1+r_2$.
\end{lemma}

\begin{proof}
Let $\phi:U_K\rightarrow \mathbb R^{r_1+r_2} $ be the group homomorphism
$$
u\mapsto (c_1\log|u^{(1)}|, \ldots, 
c_{r_1+r_2}\log|u^{(r_1+r_2)}|),
$$ 
where $c_l$ is given by \eqref{eq:cl}. Then 
$\phi(U_K) $ forms a full $(r_1+r_2-1)$-dimensional lattice in the hyperplane 
$$
H=\{{\bf e}\in \mathbb R^{r_1+r_2}: e_1+\cdots +e_{r_1+r_2}=0\}.
$$ 
It follows that for any ${\bf e}\in H$, there exists a unit $u\in U_K$ 
such that $|\phi(u)-{\bf e}|\ll 1$.

Now let $a_l=
c_l\log |v^{(l)}|$ and 
 $b_l=
c_l(\log |\nm (v)|)/d$, for $1\leq l\leq r_1+r_2$.
Then it is clear that 
${\bf e}=\a-\b \in H$.  Hence we can find 
 $w\in U_K$ such that $|\phi(w)-{\bf e} |\ll 1$, which implies that
$|\a-\phi(w)-\b|\ll 1$. The lemma follows on taking 
$u=w^{-1}$.
\end{proof}

\subsection{The Dedekind zeta function}

In this section we discuss the 
Dedekind zeta function associated to $K$. 
All of the facts that we record may be found in the work of Landau \cite{landau}, for
example.
The Dedekind zeta function is defined to be 
$$
\zeta_K(s)=\sum_{(0)\neq\fa\subseteq \fo} (\n \fa)^{-s},
$$
for any $s=\sigma+it \in
\CC$ with $\sigma=\Re(s)>1$. The zeta function admits a meromorphic
continuation to the entire complex plane with a simple pole at $s=1$. 
Let
$$
A=\frac{\sqrt{D_K}}{2^{r_2}\pi^{d/2}}.
$$
 The functional equation
for the Dedekind zeta function may be written $\Phi(s)=\Phi(1-s)$,
with
$$
\Phi(s)=A^s \Gamma\left(\frac{s}{2}\right)^{r_1}\Gamma(s)^{r_2} \zeta_K(s).
$$
One recalls that $\Gamma(s)$ has simple poles at each non-positive
integer. 
By the class number
formula we have 
\begin{equation}\label{eq:Delta}
\Delta_{K}=
\Res_{s=1} \zeta_K(s)=
\frac{h_{K}2^{r_1} (2\pi)^{r_2}R_K}{w_{K} \sqrt{D_{K}}},
\end{equation}
where  $h_{K}$ is the class number of $K$, $w_{K}$ is the number of  
roots of unity in $K$ and $R_K$ is the regulator. 
It follows from the functional equation that $\zeta_K(s)$ has a zero
of order $r_1+r_2-1$ at 
$s=0$,  zeros of order $r_1+r_2$ at all negative even integers 
and 
 zeros of order $r_2$ at all negative odd  integers.
Recalling  that $\Gamma(1/2)=\sqrt{\pi}$ it easily follows that
\begin{equation}\label{eq:Psi}
\Upsilon_{K}=
\Res_{s=0} \frac{\zeta_K(s)}{s^{r_1+r_2}}=
-
\frac{\Delta_K\sqrt{D_K}}{2^{r_1}(2\pi)^{r_2}}.
\end{equation}

We will also require some information about the order of magnitude of $\zeta_K(s)$.
The functional equation yields $\zeta_K(s)=\gamma_K(s)\zeta_K(1-s)$, with 
$$
\gamma_K(s)=
A^{1-2s} \left(\frac{\Gamma(\frac{1-s}{2})}{\Gamma(\frac{s}{2})}
\right)^{r_1}
\left(\frac{\Gamma(1-s)}{\Gamma(s)}\right)^{r_2}.
$$
In any fixed strip $\sigma_1\leq \sigma \leq \sigma_2$, an application of 
Stirling's formula yields the existence of 
$\Lambda_\sigma\in \CC$ and $\Lambda\in \RR$, depending on $K$, 
such that 
 \begin{equation}\label{eq:stirling}
\gamma_K(s)= \Lambda_\sigma t^{(1/2-\sigma)d} e^{-\pi it(\log t-\Lambda)} 
\left(1+O\left(\frac{1}{t}\right)\right),
\end{equation}
as $t\rightarrow \infty$. This is established in \cite[Satz 166]{landau}. Here 
$\Lambda=\log 2\pi +1 -\frac{1}{d}\log D_K$ is an absolute constant but $\Lambda_\sigma$ depends
on $\sigma$ and satisfies 
$|\Lambda_\sigma|= D_K^{1/2-\sigma} (2\pi)^{d(\sigma-1/2)}$. 
We have  
\begin{equation}\label{eq:mu} 
\limsup_{t\rightarrow \pm \infty} \frac{\log |\zeta_K(\sigma+it)|}{\log |t|}
\leq \begin{cases}
0, & \mbox{if $\sigma>1$,}\\
(1-\sigma)d/2, & \mbox{if $0\leq \sigma\leq 1$,}\\
(1/2-\sigma)d, & \mbox{if $\sigma< 0$.}
\end{cases}
\end{equation}
Here the first inequality is obvious and the final inequality follows from the functional equation
for $\zeta_K(s)$ and \eqref{eq:stirling}. The middle inequality is a consequence of convexity.

\subsection{Smooth weight functions}\label{s:weights}

Let $K$ be a number field of degree $d$, as previously, and let $V$ be
the associated $\RR$-algebra $\RR^{r_1}\times \CC^{r_2}$. 
When dealing with functions on $V$ or $V^n$ 
it will occasionally be convenient to work with alternative coordinates 
 $v=(\v,\y+i\z)$ on $V$, with 
$$
\v=(v_1\ldots,v_{r_1}), \quad
\y=(y_1\ldots,y_{r_2}), \quad
\z=(z_1\ldots,z_{r_2}).
$$
We will then employ the volume form
$$
\d v = \d v_1 \cdots \d v_{r_1} \d y_1  \cdots \d y_{r_2}\d z_1\cdots \d z_{r_2},
$$
on $V$, 
where $\d v_i$ is the standard Lebesgue measure on $\RR$ for
$1\leq i \leq r_1$ and $\d y_j, \d z_j$ are, respectively, the
standard Lebesgue measures on $\Re(\CC)$ and $\Im(\CC)$, for $1\leq
j\leq r_2$.  The corresponding volume form on $V^n$ will be denoted by
$\d \v$ or $\d \x$.

Our work will make prevalent use of smooth weight functions on $V$,
and more generally on $V^n$, for integer $n\geq 1$.
For us a smooth weight function on $V^n$ is any infinitely
differentiable function $w: V^n \rightarrow \CC$ which has compact support.
The latter is equivalent
to the existence of $A>0$  such that $w$ is supported on the hypercube $[-A,A]^{dn}$.  Let $n=1$.
For any $\bbe=(\beta_1,\ldots,\beta_{d})\in \mathbb{Z}_{\geq 0}^{d}$ and any smooth weight function $w$ on $V$, we will use the notation 
$$
\partial^{\bbe}
w(v)=\prod_{i=1}^{r_1}\partial_{v_i}^{\beta_i}\prod_{j=1}^{r_2}\partial^{\beta_{j+r_1}}_{y_j}
\partial^{\beta_{j+r_1+r_2}}_{z_j}w(v).
$$
We will denote the ``degree'' of $\bbe$ by $|\bbe|=\beta_1+\cdots +\beta_d$.
When $n\geq 1$ is arbitrary, 
an analogous definition of $\partial^{\bbe} w(\x)$ will be used for any smooth weight $w$ on $V^n$, for 
$\bbe\in \mathbb{Z}_{\geq 0}^{dn}$.
For a smooth weight $w$ on $V^n$, and any $N\geq 0$, we let
\begin{equation}\label{eq:sobolev}
\lambda_w^N=\sup_{\substack{\x\in V^n\\ |\bbe|\leq N}} \left| \partial^{\bbe} w(\x)\right|.
\end{equation}
We will henceforth write $\mathcal{W}_n(V)$ for the set of smooth
weight functions $w$ on $V^n$ for which 
$\lambda_w^N$ is bounded by an absolute constant $A_N>0$, for each integer $N\geq 0$.
We will write $\mathcal{W}_n^+(V)$ for the subset of $w \in
\mathcal{W}_n(V)$ which take values on non-negative real numbers
only. 
In what follows, unless explicitly indicated otherwise,  we will allow the implied constant in any estimate involving a weight $w\in
\mathcal{W}_n(V)$ to depend implicitly on $A$ and $A_{N}$.

\subsection{Additive characters over $K$}\label{s:characters}

Given any non-zero 
integral ideal $\fb$ of $K$, an additive character modulo $\fb$
is defined to be a non-zero function $\sigma$ on $\fo/\fb$ which satisfies
$$
\sigma(\alpha_{1}+\alpha_{2})=
\sigma(\alpha_{1})
\sigma(\alpha_{2}),
$$
for any $\alpha_{1},\alpha_{2}\in \fo$. Such a character is said to be primitive 
if it is not a character modulo $\fc$ for any ideal $\fc\mid \fb$, with $\fc\neq \fb$.
We will make use of the basic orthogonality relation
\begin{equation}\label{eq:orthogonal}
\sum_{\sigma \bmod{\fb}} \sigma(\alpha) = 
\begin{cases}
\n \fb, & \mbox{if $\fb\mid (\alpha)$,}\\
0, & \mbox{otherwise,}
\end{cases}
\end{equation} 
for any non-zero $\alpha\in \fo$,
where the notation 
$\sum_{\sigma \bmod{\fb}}$ means that the sum is taken over
 additive characters modulo $\fb$.
In fact 
there is an isomorphism between the additive characters modulo $\fb$ and the residue classes modulo
$\fb$
and so the number of distinct characters is $\n \fb$. 
In this isomorphism primitive characters correspond to residue classes that are relatively prime to
$\fb$. 
Hence there are $\phi(\fb)$ 
distinct primitive characters modulo  $\fb$.  It is easy to see that if $\sigma_{0}$ is a fixed
primitive character modulo $\fb$ then, as $\beta$ runs through elements of $\fo/\fb$ (respectively,
through elements of $(\fo/\fb)^*$), the functions $\sigma_{0}(\beta \cdot)$ give
all the characters 
(respectively, primitive characters)
modulo $\fb$ exactly once.

For a given integral ideal $\fb$ we now proceed to construct an explicit non-trivial primitive
character modulo $\fb$. For this we will need some preliminary algebraic facts to hand.
Recall the notation $\fd$ and $\fa_\gamma$, for the different ideal and denominator ideal,
respectively. 

\begin{lemma}\label{lem:alg1}
Let $\ve>0$ and let $\fb,\fc$ be integral ideals.
Then we have the following:
\begin{enumerate}
\item[(i)]
there exists $\alpha\in \fb$ such that $\ord_{\fp}(\alpha)=\ord_\fp(\fb)$ for every prime ideal
$\fp\mid \fc$;
\item[(ii)]
there exists $\alpha\in \fb$ and an unramified prime ideal $\fp$ coprime to $\fb$, 
with $\n\fp \ll (\n \fb)^{\ve}$, 
 such that $(\alpha)=\fb\fp$.
\end{enumerate}
\end{lemma}

\begin{proof}
Part (i) is standard. 
For part (ii) we
note that $\fb$ has $\omega(\fb)=O(\log \n \fb)$ prime ideal divisors. However the number of
distinct prime ideals with norm at most $(\n \fb)^\ve$ is $\gg (\n \fb)^{\ve/2}$. Using the
equidistribution of prime ideals in ideal classes we therefore deduce that there exists an
unramified prime ideal $\fp$ such that $\n \fp \ll (\n \fb)^\ve$, with $\fp\nmid \fb$ and 
$\fb\fp$ principal. Any generator $\alpha$ of $\fb\fp$ then satisfies the properties claimed. 
\end{proof}

The  $\alpha$ constructed in part (ii) satisfies 
$\ord_{\fq}(\alpha)=\ord_\fq(\fb)$ for every prime ideal
$\fq\mid \fb$, with  $|\n_{K/\QQ}(\alpha)|\ll (\n \fb)^{1+\ve}$.
 Thus part (ii) is a refinement of part (i) in 
the special case $\fc=\fb$, 
with additional control over the size of
the norm of $\alpha$.
We now have everything in place to construct our non-trivial character modulo an integral ideal
$\fb.$ Consider the integral ideal
$$
\fc=\fb\fd. 
$$
By Lemma \ref{lem:alg1} (ii) we can find an integer $\alpha\in \fc$ and an unramified prime
ideal $\fp_1$ coprime to $\fb$,  with $\n\fp_1 \ll (\n \fb)^{\ve}$, such that
$(\alpha)=\fc\fp_1$.
Applying  Lemma \ref{lem:alg1} (i) 
we see that there exists $\nu\in \fo$ such that $\fp_1 \mid (\nu)$
but $(\nu)$  and $\fc$ are coprime.
It follows that $\fc=\fa_\gamma$ with $\gamma=\nu/\alpha$.
Indeed, $\beta\in \fa_\gamma$ if and only if $(\beta \nu)\subseteq (\alpha)=\fc \fp_1$, which is 
if and only if $\beta\in \fc$.
We claim that 
$$
\sigma_0(\cdot) = \e(\gamma \cdot ) 
$$
defines a non-trivial primitive character modulo $\fb$.

To check the claim we note that 
$\sigma_0$ is a trivial character if and only if $\gamma\in \hat \fo$. But this holds if and only
if $\fa_\gamma\supseteq  \fd$, which is so if and only if $\fc=\fa_\gamma \mid \fd$. This is clearly
impossible. Suppose now that $x,z\in \fo$, with $\fb\mid (z)$. Then it is clear that 
$\fa_\gamma \mid (z)\fd$, whence 
$\gamma z\in \hat \fo$. Thus $\sigma_0(x+z)=\sigma_0(x)$ and so it follows that $\sigma_0$ is a non-trivial
character modulo $\fb$. Lastly, we need to check the primitivity of the character. Let $\fb_1\mid \fb$
be any proper ideal divisor. Then the ideal $\fb_1\gamma$ is not contained in $\hat \fo$ and so there exists  
$z\in \fb_1$ such that $\gamma z\not \in \hat \fo$.
 This implies
that $\sigma_0$ is not a character modulo any ideal $\fb_1\mid \fb$, with $\fb_1\neq \fb$,
so that it is in fact  a primitive character modulo $\fb$.

We may summarise our investigation in the following result.

\begin{lemma}\label{lem:orthogonal}
Let $\ve>0$ and let $\fb$ be an integral ideal. Then there exists $\gamma\in K$, with
$\gamma=\nu/\alpha$
 for 
$\alpha\in \fb\fd$ and $\nu\in \fo$ 
such that $(\nu)$ is coprime to $\fb\fd$, 
together with  a prime ideal $\fp_1\mid (\nu)$ satisfying
$\n \fp_1 \ll (\n \fb)^{\ve}$ and $(\alpha)=\fb\fd\fp_1$,  such that 
$\e(\gamma\cdot)$ defines a non-trivial primitive additive character modulo $\fb$. In particular, we have 
$$
\starsum_{\sigma\bmod{\fb}
}\sigma(x) = 
\sum_{a\in (\fo/\fb)^*} \e(a \gamma x),
$$
for any $x\in \fo$.
\end{lemma}

Finally, we need to discuss how additive  characters modulo $\fb$ can be extended to arbitrary  integral ideals $\fa\subseteq \fo$, as in the statement of Theorem \ref{t:delta}. 
By part (i) of Lemma  \ref{lem:alg1} there exists 
$\alpha \in \fo$ such that $\ord_\fp(\alpha)=\ord_\fp(\fa)$ for all
$\fp \mid \fa\fb\fd$.  Given any additive character modulo $\fb$ we will extend it to $\fa$ by setting
$
\sigma(\fa)=\sigma(\alpha).
$
On noting that  $\fb\mid (\alpha) $
if and only if $\fb\mid \fa $,  we conclude from
\eqref{eq:orthogonal} that  
\begin{equation}\label{eq:orthogonal'}
\sum_{\sigma \bmod{\fb}} \sigma(\fa) = 
\begin{cases}
\n \fb, & \mbox{if $\fb\mid \fa$,}\\
0, & \mbox{otherwise.}
\end{cases}
\end{equation} 
We will need to know that the sum  
$$
\starsum_{\sigma\bmod{\fb}
}\sigma(\fa),
$$
over primitive characters modulo $\fb$, 
does not depend on the choice of $\alpha$.   This is achieved in the following result. 

\begin{lemma}
\label{delta1}
Let $\fa$ be an integral ideal. 
Let $\alpha_1, \alpha_2\in \fo$  such  that
$\ord_\fp(\alpha_i)=\ord_\fp(\fa)$ for $i=1,2$ and  all
$\fp\mid \fa\fb\fd$. Then
$$\starsum_{\sigma\bmod{\fb}}\sigma(\alpha_1)=\starsum_{\sigma\bmod{\fb}
}\sigma(\alpha_2).
$$
\end{lemma}

\begin{proof}
Let $\fa_1,  \fa_2$ be coprime integral ideals, which are also
coprime to $\fa\fb\fd$, such that 
$(\alpha_1^{-1}\alpha_2)=\fa_1^{-1}\fa_2$.
Since $\fa_1^{-1}\fa_2$ is principal, the ideals $\fa_1$ and $\fa_2$ are in
the same ideal class. Let $\fp_1$ be a prime ideal, coprime to $\fa\fb\fd $, such that
$\fa_1\fp_1$ and $\fa_2\fp_1$ are principal ideals. Thus we can choose integers $\beta_1$ and $\beta_2$
such that $\fa_1\fp_1=(\beta_1)$ and $\fa_2\fp_1=(\beta_2)$, with  
$\beta_1\alpha_2=\beta_2\alpha_1$.  
In particular $\beta_1$ and $\beta_2$ are coprime to $\fb\fd$ and it follows that 
\begin{align*}
 \starsum_{\sigma\bmod{\fb}}\sigma(\alpha_1)
 &=\starsum_{\sigma\bmod{\fb}}\sigma(
\beta_2\alpha_1)\\
&=\starsum_{\sigma\bmod{\fb}}\sigma(\beta_1\alpha_2)\\
&=\starsum_{\sigma\bmod{\fb}}\sigma(\alpha_2),
\end{align*}
as required. 
\end{proof}

It is often  convenient to restrict a sum over additive characters to a sum over primitive characters,
using the identity
\begin{equation}\label{eq:south}
\sum_{\sigma\bmod{\mathfrak{c}}}\sigma(\fa)=\sum_{\mathfrak{b}\mid \fc}
\sideset{}{^{*}}\sum_{\sigma\bmod{\mathfrak{b}}}\sigma(\fa).
\end{equation}
The previous lemma ensures that this is well-defined.

\subsection{Counting rational points on algebraic varieties}

Let $G\in \fo[X_1,\dots, X_n]$
be a 
 homogeneous polynomial, 
 which is absolutely irreducible over $K$ and has degree  $\geq 2$. 
We will need an estimate for the number 
of $\x\in \fo^n$ for which $G(\x)=0$, subject to certain constraints. 

\begin{lemma}\label{lem:cohen}
Let $B=\bigoplus_l B^{(l)} \in V$,  with $B^{(l)}>0$
and  $\nm(B)\geq 2$. Let $\ve>0$,  let $\fc$ be an integral ideal and 
let  $\a\in \fo^n$. Then we have 
\[
\#\left\{\x\in \fo^n:  
\begin{array}{l}
| \x^{(l)}|\leq B^{(l)}, \\
G(\x)=0, ~ \x \equiv \a \bmod{\fc}
\end{array}
\right\}
\ll  (\nm(B))^\ve \left(1+\frac{\nm(B)}{\n\fc}\right)^{n-3/2}.
\]
\end{lemma}

The implied constant in this estimate depends at most on $G$ and the choice of $\ve$.
Lemma~\ref{lem:cohen} is a generalisation of 
\cite[Lemma 15]{hb-10} to the number field setting. One expects that 
one should be able to replace the exponent $n-3/2$ by $n-2$, whereas in fact any exponent less than  $n-1$ suffices to obtain Theorem \ref{t:10}.

\begin{proof}[Proof of Lemma \ref{lem:cohen}]
We denote by $N(B;\fc,\a)$ the quantity that is to be estimated.
For the proof we may assume without loss of generality that $\n\fc\leq \nm(B)/2$, for otherwise the result
follows from the trivial bound $N(B;\fc,\a)\ll(1+\nm(B)/\n\fc)^{n}$.
Let us define $N(B;\fc)=\max_{\a\in \fo^n}N(B;\fc,\a)$ and let $\delta>0$.

Applying Lemma \ref{lem:alg1}(ii) we see that there exists 
$\gamma\in \fo$ and an unramified prime ideal 
$\fp$ coprime to $\fc$,  such that
$(\gamma)=\fp\fc$ and $\n \fp\ll  (\n \fc)^{\delta}.$ Let us set $\fg=(\gamma)$.
Then it follows that 
\begin{equation}\label{eq:apple}
N(B;\fc,\a)=\sum_{\b\in (\fc/\fg)^n} N(B;\fg,\a+\b)\ll (\n\fc)^{\delta n}N(B;\fg),
\end{equation}
since $\#\fc/\fg=[\fc:\fg]=[\fo:\fg]/[\fo:\fc]=\n \fp$.
For any $u\in U_K$, where 
$U_K$  denotes the group of units of $\fo$,
 we let $\tilde u=\bigoplus_l |u^{(l)}|$. 
Then we have 
\begin{align*}
N(B;\fg)
&\leq \max_{\a\in \fo}
\#\left\{u\x\in \fo^n:  
\begin{array}{l}
| u^{(l)}\x^{(l)}|\leq |u^{(l)}|B^{(l)}, \\
G(u\x)=0, ~ u\x \equiv u\a \bmod{\fg}
\end{array}
\right\}\\
&\leq \max_{\a\in \fo}
\#\left\{\y\in \fo^n:  
\begin{array}{l}
| \y^{(l)}|\leq (\tilde u B)^{(l)}, \\
G(\y)=0, ~ \y \equiv \a \bmod{\fg}
\end{array}
\right\}\\
&=N(\tilde u B;\fg).
\end{align*}
According to Lemma \ref{lem:unit} it therefore suffices to 
estimate the quantity $N(B;\fg)$, for $B$ satisfying $
B^{(l)}\ll \nm(B)^{1/d}$.
Likewise, a further application of Lemma \ref{lem:unit} allows us to assume that 
$\fg=(\gamma)$, with 
$(\n\fg)^{1/d}\ll |\gamma^{(l)}|\ll
(\n\fg)^{1/d}$.
Now let $\x$ be any vector in $V^n$. Then there exists  $\y\in \fo^n$ such that 
$\langle \gamma^{-1}\x-\y\rangle \ll 1$.
This implies that $\langle \x-\gamma \y\rangle \ll (\n\fg)^{1/d}$, whence given any 
$\x\in V^n$ we can find a vector $\y\in  \fg^n=\gamma\fo^n$ 
such that $\langle\x-\y\rangle\ll (\n\fg)^{1/d}$.
Therefore, in our analysis of 
$N(B;\fg,\a)$ for given $\a\in \fo^n$, it suffices to assume that 
$\langle\a\rangle \ll (\n\fg)^{1/d} $.

Let $Y\subseteq \AA_K^n$ be the hypersurface 
$G(\a+\gamma\x)=0$. This is clearly absolutely irreducible of degree at least $2$.
Given any $H\geq 2$, it follows from an application of the large sieve
(see Serre \cite[Chap.~13]{serre})
that 
$$
\#\{\x\in \fo^n\cap Y:\langle\x\rangle\leq H\}\ll 
 H^{(n-3/2)d}\log H.
$$
 Furthermore, an inspection of the proof reveals that the implied constant is independent 
 of $\a$ and $\gamma$.
It now follows that
\begin{align*}
 N(B;\fg,\a)
 &=
 \#\left\{\x\in \fo^n:  
\begin{array}{l}
| (\a+\gamma \x)^{(l)}|\leq B^{(l)}, \\
G(\a+\gamma \x)=0
\end{array}
\right\}\\
 &\leq
 \#\left\{\x\in \fo^n\cap Y:  
| \x^{(l)}|\ll H
\right\}\\
&\ll H^{(n-3/2)d}\log \nm(B),
\end{align*}
where $H=(1+\nm(B)/\n\fg)^{1/d}$. 
Inserting this into \eqref{eq:apple} and taking any $\delta <\ve/n$ therefore leads to the conclusion of the lemma.
\end{proof}

\subsection{Poisson summation over $K$}

We will need a version of the Poisson summation formula for number
fields.  This is provided for us by  
the work of Friedman and Skoruppa  \cite{f-s}, in which we take the
base field to be $\QQ$ (and so $m=1$).  
On $\RR_{>0}$ we 
define the function
\begin{align*}
k_{r_1,r_2}(t)
&=
\frac{1}{2\pi i} \int_{c-i\infty}^{c+i\infty}
t^{-z} 
\left(
\frac{\pi^{-z} \Gamma(\frac{z}{2})}{\Gamma(\frac{1-z}{2})}
\right)^{r_1}
\left(
\frac{(2\pi)^{-2z}\Gamma(z)}{\Gamma(1-z)}
\right)^{r_2}
\d z\\
&=
\frac{1}{2\pi i} \int_{c-i\infty}^{c+i\infty}
t^{-z} g(z)^{r_1} h(z)^{r_2} \d z,
\end{align*}
say,   for any $0<c<1/6$.
Then $k_{r_1,r_2}$ is a Mellin convolution 
$k_{1,0}^{*r_1} *k_{0,1}^{*r_2}$, with 
$$
k_{1,0}(t)=\frac{2}{\sqrt{\pi}} \cos (2\pi t), \quad 
k_{0,1}(t)=J_{0}(4\pi \sqrt{t}).
$$
Let $f\in \cW_1(\RR)$ and let $\fa$ be a fractional ideal of
$K$. Then  the version of Poisson summation that we need takes the
form 
\begin{equation}\label{eq:poisson}
\sum_{\substack{\alpha\in \fa/U_K\\ \alpha\neq 0}}
f(|\n_{K/\QQ}\alpha|) = \frac{\Delta_{K}}{h_{K}\n \fa} \int_{0}^{\infty} f(t) \d t
+
\frac{2^{r_2}\pi^{d/2}}{\sqrt{D_{K}} \n \fa}
\sum_{\substack{\beta\in \hat \fa/U_K\\ \beta\neq 0}}
\tilde f(|\n_{K/\QQ} \beta|),
\end{equation}
where $\Delta_{K}$ is given  by \eqref{eq:Delta} and 
$$
\tilde f(y)= \int_{0}^{\infty} f(t) k_{r_1,r_2}(ty) \d t,
$$
as a function on $\RR_{>0}$.
In fact we have 
\begin{equation}\label{eq:decay}
\tilde f(y) \ll_{N} y^{-N},
\end{equation}
for any $N\geq 0$. This follows on noting that 
\begin{align*}
\tilde f(y)&=
\frac{1}{2\pi i} 
\int_{c-i\infty}^{c+i\infty}  y^{-z} g(z)^{r_1} h(z)^{r_2} 
 \int_{0}^{\infty} f(t) t^{-z}\d t \d z\\
 &=\frac{1}{2\pi i} 
\int_{1-c-i\infty}^{1-c+i\infty}  F(w) y^{w-1} g(1-w)^{r_1} h(1-w)^{r_2} \d w,
\end{align*}
where $F$ is the Mellin transform of $f$. Since $f$ is smooth and compactly
supported, it follows that $F$ is entire and of rapid decay, enough to
counter the polynomial growth of the functions $g$ and $h$. The
estimate \eqref{eq:decay} follows on  
shifting the contour sufficiently far to the left.

Using \eqref{eq:poisson} we may deduce a corresponding version in
which the sum over elements is replaced by a sum over ideals. In fact,
we will often be  led to consider sums of the form
$$
\sum_{\mathfrak{b}\neq (0)} f(R\n \mathfrak{b}),
$$
for a parameter $R>0$, with $f\in \cW_1(\RR)$ and the sum being taken over
integral ideals.  To handle this sum 
we let $[\fc_{1}],\dots,[\fc_{h_{K}}]$ 
denote distinct cosets for the class group $C(K)$.
For each $1\leq j\leq h_{K}$, 
we have a bijection between the set $S_{j}$ of
integral ideals in $[\fc_{j}]$ and the elements of the ideal
$\fc_{j}^{-1}$ modulo the action of 
the unit group  $U_K$.
The explicit bijection is  $\psi_{j}: S_{j}\rightarrow
\fc_{j}^{-1}/U_K$, given by  $\psi_{j}(\fb)=\fc_{j}^{-1}\fb$.
Using our decomposition into
ideal classes we may therefore write
$$
\sum_{\mathfrak{b}\neq (0)} f(R\n \mathfrak{b})
=\sum_{j=1}^{h_{K}}
\sum_{\substack{\alpha\in \fc_j^{-1}/U_K\\ \alpha\neq 0}} f(R\n (\alpha\fc_j)). 
$$
Let $f_j(t)=f(Rt\n \mathfrak c_j)$.  
We analyse the inner sum using the  Poisson summation formula 
in the form \eqref{eq:poisson}. This gives
 \begin{align*}
\sum_{\substack{\alpha\in \fc_j^{-1}/U_K
\\ \alpha\neq 0}
} 
f_{j}(|\n_{K/\QQ}(\alpha)|)
&=
 \frac{\Delta_{K} \n \fc_{j}}{h_{K}} \int_{0}^{\infty} f_{j}(t) \d t
+
\frac{2^{r_2}\pi^{d/2} \n \fc_{j}}{\sqrt{D_{K}} }
\sum_{\substack{\beta\in \hat \fc_{j}^{-1}/U_K\\ \beta\neq 0}}
\tilde f_{j}(|\n_{K/\QQ} (\beta)|).
\end{align*}
A change of variables reveals that
$$
\int_{0}^{\infty} f_{j}(t) \d t
=\frac{1}{R\n \fc_{j}}
\int_{0}^{\infty} 
f(t)\d t.
$$
Moreover, it follows from \eqref{eq:decay} that 
\begin{align*}
\tilde f_{j}(y) 
&= \int_{0}^{\infty} 
f( Rt\n \mathfrak c_j)
 k_{r_1,r_2}(ty) \d t\\
 &= \frac{1}{R\n \fc_{j}}
\tilde f\left(\frac{y}{R\n \fc_{j}}\right) \\
&\ll_{N} R^{-1}\left(\frac{y}{R}\right)^{-N},
\end{align*}
for any  $N\geq 0$. Reintroducing the sum over $j$, we
have therefore established the following result.

\begin{lemma}\label{lem:poisson}
Let $f\in \cW_1(\RR)$ and let $R>0$. Then we have 
$$
\sum_{\mathfrak{b}\neq (0)} f(R\n \mathfrak{b})=
\frac{\Delta_{K}}{R}
\int_{0}^{\infty} 
f(t)\d t +O_N\left(R^N  \right),
$$
for any $N>0$.
\end{lemma}

\section{The smooth $\delta$-function over $K$}\label{s:delta}

In this section we establish Theorem \ref{t:delta} and some further
basic properties of the function $h(x,y)$. 
Let $w$ be any infinitely differentiable bounded  
non-negative function on $\mathbb R$ which is
supported in the interval $[1/2,1]$ and satisfies
$$
\int_{-\infty}^\infty w(t)\d t=1.
$$
In particular $w\in \cW_1^+(\RR)$.
For $Q\geq 1$ define 
$$
c_Q^{-1}=\Delta_K^{-1}Q^{-d}\sum_{\mathfrak{c}}w(Q^{-d}\n \mathfrak{c}  ) ,
$$ 
where $\Delta_{K}$ is given by \eqref{eq:Delta} and the sum is over integral ideals $\fc\subseteq
\fo$. 
Here, as throughout this section, the support of $w$ restricts the sum to non-zero ideals.
For any non-zero ideal $\fa\subseteq \fo$, we have
$$
\sum_{\mathfrak{c}  \mid \fa}\left\{w\left(\frac{\n \mathfrak{c}}{Q^{d}}  \right)-
w\left(\frac{\n \fa}{Q^{d}\n \mathfrak{c}}\right)\right\}=0,
$$
where the sum is over  integral ideals $\mathfrak{c}  $ dividing $\fa$.
On the other hand, when $\fa=(0)$, we have 
$$
\sum_{\mathfrak{c}  \mid \fa}\left\{w\left(\frac{\n \mathfrak{c}}{Q^{d}}  \right)-
w\left(\frac{\n \fa}{Q^{d}\n \mathfrak{c}}\right)\right\}=
\sum_{\mathfrak{c}}w(Q^{-d}\n \mathfrak{c}  ). 
$$ 
In this way we deduce that
\begin{equation*}
\delta_K(\fa)=
c_Q
\Delta_K^{-1}Q^{-d}
\sum_{\mathfrak{c}  \mid \fa}\left\{w\left(\frac{\n \mathfrak{c}}{Q^{d}}  \right)-
w\left(\frac{\n\fa}{Q^{d}\n \mathfrak{c}}\right)\right\}.
\end{equation*}
Using \eqref{eq:orthogonal'} to detect the divisibility condition, we may therefore write
\begin{equation}
\label{alph1}
\delta_K(\fa)=
c_Q
\Delta_K^{-1} Q^{-d}
\sum_{\mathfrak{c}}\frac{1}{\n \fc}
\sum_{\sigma\bmod{\mathfrak{c}}}\sigma(\fa)
\left\{w\left(\frac{\n \mathfrak{c}}{Q^{d}}  \right)-
w\left(\frac{\n \fa}{Q^{d}\n \mathfrak{c}}\right)\right\}.
\end{equation}
Inserting \eqref{eq:south} into \eqref{alph1} and re-ordering the summation we arrive at the expression for 
$\delta_{K}(\fa)$ in Theorem \ref{t:delta}, with 
\begin{equation}\label{eq:h}
h(x,y)=\frac{1}{\Delta_K}\sum_{\mathfrak{z}}\frac{1}{x\n
  \mathfrak{z}}\left\{w(x\n \mathfrak{z})-w\left(\frac{|y|}{x\n
      \mathfrak{z}}\right)\right\}. 
\end{equation} 
Since $w$ is bounded and supported in $[1/2,1]$ we see that 
$h(x,y)\neq 0$ only if $x\leq \max\{1,2|y|\} $ and, furthermore,  
$$
h(x,y)\ll \frac{1}{x} \left( \sum_{(2x)^{-1}\leq \n \fz \leq x^{-1}}
  \frac{1}{\n \fz}
+
\sum_{|y|/x \leq \n \fz \leq 2|y|/x}
  \frac{1}{\n \fz}\right) \ll \frac{1}{x}.
$$
Finally, $h(x,y)$ is easily seen to be infinitely differentiable  on
$(0,\infty)\times \RR$ since the individual summands are infinitely
differentiable. 

In order to complete the proof of Theorem \ref{t:delta}, it 
remains to show that 
\begin{equation}
c_Q=1+O_N(Q^{-N} )\label{cq}
\end{equation} 
for any  $N>0$.
For this we apply Lemma \ref{lem:poisson} with $R=Q^{-d}$ and $f=w$, concluding that
$$
\sum_{\mathfrak{c}}w(Q^{-d}\n \mathfrak{c}  )
=
\Delta_{K}Q^{d}
\int_{0}^{\infty} 
w(t)\d t +O_N\left(Q^{-N}  \right)
$$
for any $N>0$. The first term here is $\Delta_{K}Q^{d}$, whence
\begin{align*}
c_{Q}^{-1}
&=
\Delta_{K}^{-1}Q^{-d}
\sum_{\mathfrak{c}}w(Q^{-d}\n \mathfrak{c}  )
=1+O_{N}(Q^{-N}),
\end{align*}
as claimed.

We now turn to a detailed analysis of the function $h(x,y)$ in \eqref{eq:h}. 
If $a_{m}$ are the coefficients appearing in the Dedekind zeta
function $\zeta_{K}(s)$, then we may write
$$
h(x,y)=\frac{1}{\Delta_K}\sum_{m=1}^\infty
\frac{a_{m}}{xm}\left\{w(xm)-w\left(\frac{|y|}{xm}\right)\right\}. 
$$
Our task is to achieve analogues of the corresponding 
facts established by Heath-Brown \cite[\S 4]{H} concerning $h(x,y)$ in the case $K=\QQ$, which
corresponds to taking $a_{m}=1$ for all $m$. 
However, rather than the Euler--Maclaurin formula, which is used
extensively by Heath-Brown, we will use the Poisson summation formula
in the form Lemma \ref{lem:poisson}.   The basic structure of the
proofs will nonetheless remain  similar. We begin with the
following result. 

\begin{lemma}\label{lem:vanish}
The function $h(x, y)$ vanishes when $x\geq 1$ and $|y|\leq x/2$. When
$x\leq 1$ and $|y|\leq x/2$ the function $h(x, y)$ is constant with
respect to $y$, 
taking the value 
$$
\frac{1}{\Delta_K}\sum_{\mathfrak{z}}\frac{1}{x\n \mathfrak{z}}w(x\n \mathfrak{z}).
$$
\end{lemma}

\begin{proof}
This is obvious  on
recalling that $w$ has  support $[1/2,1]$.  
\end{proof}

Our next task is to show that the sum involved in $h(x,y)$ nearly cancels if
$x=o(\min(1,|y|))$, for which we will closely follow
the argument used in \cite[Lemma 4]{H}.
We will require some preliminary
estimates for $h(x,y)$ and its partial derivatives with respect to $x$
and $y$.  
In view of \eqref{eq:h} we may write $h(x,y)=h_1(x,y)-h_2(x,y)$, 
with 
\begin{equation}\label{eq:hi}
h_k(x,y)=\frac{1}{\Delta_K}\sum_{\mathfrak{z}} f_k (\n \mathfrak{z})
\end{equation}
for $k=1,2$, where
 \begin{align*}
 f_1(t)
=
 \frac{1}{xt }
 w(xt),\quad
 f_2(t)=
 \frac{1}{xt }
  w\left(\frac{|y|}{xt}\right).
 \end{align*}
Let $i,j\in \ZZ_{\geq 0}$.
A simple induction argument reveals that 
\begin{equation}
\frac{\partial^{i+j}}{\partial x^i\partial
y^j} \left( 
\frac{1}{x}w\left(\frac{y}{xn}\right)
\right)
=x^{-i-1} y^{-j}\sum_{0\leq t
\leq i}c_{i,j,t}\left(\frac{y}{xn}\right)^{j+t}w^{(j+t)}\left(\frac{y}{xn}\right)\label{eq:y derivative}
\end{equation}
and 
\begin{equation}\label{eq:x derivative}
\frac{\partial^{i}}{\partial x^i}\left(\frac{1}{x} w(x n)\right)
=\sum_{0\leq t\leq i} d_{i,t}x^{-i+t-1} n^{t}w^{(t)}(x n),
\end{equation}
for certain constants $c_{i,j,t}=O_{i,j}(1)$ and $d_{i,t}=O_i(1)$, with $c_{i,0,0}=d_{i,0}$
and $c_{0,0,0}=d_{0,0}=1$.

Suppose that 
$x\leq 1$ and $|y|\leq x/2$. 
In particular 
$n^{t}w^{(t)}(x n)\ll x^{-t}$ whenever $1/2\leq xn\leq 1$.   
Thus it follows from 
Lemma \ref{lem:vanish} and \eqref{eq:x derivative} that
\begin{equation}
  \label{eq:1}
\frac{\partial^i}{\partial x^i} h(x,y) \ll_{i} 
x^{-i-1} \sum_{(2x)^{-1}\leq \n \fz\leq x^{-1}} \frac{1}{\n \fz} \ll_i x^{-i-1}.
\end{equation}
Next suppose that  $|y|\geq x/2$. We claim that 
\begin{equation}
  \label{eq:2}
\frac{\partial^{i+j}}{\partial x^i\partial y^j} h(x,y) \ll_{i,j} x^{-i-1}|y|^{-j}.
\end{equation}
Now the terms $(x\n\fz )^{-1} w(x\n\fz )$
in the definition \eqref{eq:hi} of $h_1(x, y )$ only contribute to the partial derivative when $x\leq 1$
and $j = 0$, in which case they contribute $O_i(x^{-i-1})$ as before, which is satisfactory.
For the terms in $h_{2}(x,y)$ we apply \eqref{eq:y derivative}. Since
\[
\left(\frac{y}{xn}\right)^{j+t}w^{(j+t)}\left(\frac{y}{xn}\right)\ll_{j,t} 1,
\]
it follows that
\[
\frac{\partial^{i+j}}{\partial x^i\partial y^j} h_2(x,y)\ll_{i,j}
x^{-i-1}|y|^{-j} \sum_{x/|y|\leq\n\mathfrak{z}\leq
2x/|y|}\frac{1}{\n \mathfrak{z}}\ll_{i,j}
x^{-i-1}|y|^{-j},\]
which is satisfactory.

\begin{lemma}\label{lem:decay}
Let $i, j,N \in \ZZ_{\geq 0}$. Then we have
$$
\frac{\partial^{i+j}}{\partial x^i\partial y^{j}} h(x,y)
\ll_{i,j, N}x^{-i-j-1}\left(x^{N}+\min\{1,(x/|y|)^{N}\}\right).
$$ 
The term $x^N$ on the right can be omitted if $j\neq 0$.
\end{lemma}

\begin{proof}
If $x\geq 2|y|$ then the result is an immediate consequence of
Lemma \ref{lem:vanish} and 
\eqref{eq:1}. If $|y|\leq x\leq 2|y|$ then it follows from 
\eqref{eq:2}. Hence we may assume that $x\leq |y|$.
The case $N=0$ is trivial and so we suppose that $N\geq 1$. 

Let $i,j\in \ZZ_{\geq 0}$.  
Our argument is similar to the proof of \eqref{cq}, being based on
Lemma \ref{lem:poisson}. 
Writing $h(x,y)=h_1(x,y)-h_2(x,y)$, as in \eqref{eq:hi},
we deduce from \eqref{eq:x derivative} that 
\begin{align*}
\frac{\partial^{i}}{\partial
x^i} h_1(x,y)
&=\frac{1}{\Delta_K}\sum_{\fz}(\n\fz)^{-1}\frac{\partial^{i}}{\partial x^i}\left( \frac{1}{x}w(x\n\fz)\right)\\
&=\frac{1}{\Delta_K}x^{-i} \sum_{0\leq t\leq i} d_{i,t} \sum_{\fz}(x\n\fz)^{t-1}w^{(t)}(x\n\fz).
\end{align*}
We may assume that $x\leq 1$, else $h_1(x,y)=0$.
Let $g_{t}(z)=z^{t-1}w^{(t)}(z)$. On applying integration by parts repeatedly, we deduce that
\begin{equation}\label{eq:hallow}
\int_0^\infty  g_{t}(z)\d z=0,
\end{equation}
for $t\geq 1$.
Calling upon  Lemma \ref{lem:poisson} with $f=g_{t}$ and $R=x$, we obtain
\begin{equation}\label{eq:ween}
\begin{split}
  \frac{\partial^{i}}{\partial
x^i} h_1(x,y)
&=x^{-i-1}\sum_{0\leq t
\leq i}d_{i,t}\left( \int_0^\infty g_{t}(z)\d z+O_N(x^N)\right) \\
&=x^{-i-1}\left(d_{i,0}\int_0^\infty \frac{w(z)}{z} \d z+O_N(x^N)\right).
\end{split}
\end{equation}
Here we recall that $w$ is supported on the interval $[1/2,1]$, so that the latter integral is well-defined. 

Turning to the term $h_2(x,y)$ in \eqref{eq:hi}, we 
deduce from 
\eqref{eq:y derivative} that 
\begin{align*}
\frac{\partial^{i+j}}{\partial
x^i\partial y^{j}} h_2(x,y)
&=\frac{1}{\Delta_K}
\sum_{\fz}(\n\fz)^{-1}\frac{\partial^{i+j}}{\partial x^i\partial y^{j}}\left(\frac{1}{x}w\left(\frac{|y|}{x\n\fz}\right)\right)\\
&=\frac{1}{\Delta_K}x^{-i-1} |y|^{-j}\sum_{\fz}(\n\fz)^{-1}
\sum_{0\leq t
\leq i}c_{i,j,t}\left(\frac{|y|}{x\n \fz}\right)^{j+t}w^{(j+t)}\left(\frac{|y|}{x\n \fz}\right)\\
&=\frac{1}{\Delta_K}
x^{-i} |y|^{-j-1}
\sum_{0\leq t
\leq i}c_{i,j,t}
\sum_{\fz}
\left(\frac{|y|}{x\n \fz}\right)^{j+t+1}w^{(j+t)}\left(\frac{|y|}{x\n \fz}\right).
\end{align*}
Let $f_{j,t}(z)=z^{-j-t-1}w^{(j+t)}(1/z)$ be functions on $(0,\infty)$. By making the change of variables 
$u=1/z$, it  easily follows from \eqref{eq:hallow}  that 
\begin{align*}
 \int_0^\infty f_{j,t}(z)\d z
&=\int_0^\infty u^{j+t-1}w^{(j+t)}(u)\d u =0,
\end{align*}
when $j+t\geq 1$.
Applying Lemma \ref{lem:poisson} with $f=f_{j,t}$ and $R=x/|y|$, we get
\begin{align*}
  \frac{\partial^{i+j}}{\partial
x^i\partial y^{j}} h_2(x,y)
&=x^{-i} |y|^{-j-1}
\sum_{0\leq t
\leq i}c_{i,j,t}
\left(
\frac{|y|}{x}
\int_0^\infty f_{j,t}(z)\d z+O_N\left((|y|/x)^{-N}\right) \right)\\
&=x^{-i-1} |y|^{-j}
\left(c_{i,j,0}\int_0^\infty f_{j,0}(z)\d z+O_N\left((|y|/x)^{-N}\right) \right).
\end{align*}
When $j\neq 0$ the integral vanishes and this estimate is satisfactory for the lemma. 
On the other hand, when $j=0$ this becomes
\begin{align*}
  \frac{\partial^{i}}{\partial
x^i} h_2(x,y)
&=x^{-i-1} 
\left(c_{i,0,0}\int_0^\infty \frac{w(z)}{z}\d z+O_N\left((|y|/x)^{-N}\right) \right).
\end{align*}
Once combined with 
\eqref{eq:ween} and the fact that $c_{i,0,0}=d_{i,0}$, 
this therefore concludes the proof of the lemma.
\end{proof}

A crucial step in the proof of Lemma \ref{lem:decay} 
involved writing $h(x,y)=h_1(x,y)-h_2(x,y)$.  Here  $h_1(x,y)$
satisfies the asymptotic formula \eqref{eq:ween} with $i=0$, which was  proved under the assumption that $x\leq 1$,
but continues to hold  when $x>1$. It follows that 
\begin{equation}\label{eq:Hxy}
h(x,y)=x^{-1}H(|y|/x)+O_N(x^{N}),
\end{equation}
for any $N>0$, where $H:\RR_{>0}\rightarrow \RR$ is given by 
$$
H(v)=
\int_0^\infty \frac{w(t)}{t} \d t-
\frac{1}{\Delta_K}
\sum_{\mathfrak{z}}\frac{1}{\n \fz} w\left(\frac{v}{\n \mathfrak{z}}\right).
$$

\section{Weighted averages of $h(x,y)$}

A key ingredient in our main term analysis will be a suitable
variant of \cite[Lemma 9]{H}, showing that for small values of $x$ the
function $h(x,y)$ acts like a $\delta$-function. In point of fact we
shall be interested in the weighted average
$$
I(x)=\int_V f(v) h(x,\nm (v)) \d v,
$$
for $0<x\ll 1$ and $f\in \cW_1(V)$.
Recall from \S \ref{s:weights} that associated to each $f\in \cW_1(V)$ is a sequence of constants $A$ and $A_N$.
In a departure from our earlier conventions, in this section we will need to keep track of the dependence on the $A_N$ in any implied constant (preserving the convention that any implied constant is allowed to depend on $A$).

When $K=\QQ$ one finds that $I(x)$ is approximated by $f(0)$ to
within an error of $O_{f, N}(x^N)$, for any $N> 0$. This is achieved
through multiple applications of the Euler--Maclaurin summation
formula, an approach that is not readily adapted to the setting of
general  $K$. Instead, we will argue using   Mellin
transforms.  
Recall the definition \eqref{eq:sobolev} of $\lambda_f^N$, for any smooth weight function 
$f:V\rightarrow \CC$ with compact support.
Our goal in this section is a proof of the following
result.

\begin{lemma}\label{prop:1}
Let $f\in \cW_1(V)$ and let $N> 0$. Then we have 
$$
I(x)=
\frac{\sqrt{D_K}}{2^{r_2}}
f(0)
+
O_N\left(\lambda_f^{2d(N+1)} x^N\right).
$$
\end{lemma}

We have not attempted to obtain a dependence on $\lambda_f^M$, with $M$ minimal, since all that is required in our application is that $M$ be polynomial in $N$.
When $K=\QQ$, so that $D_K=1$ and $r_2=0$, we retrieve \cite[Lemma 9]{H}.
Since $h(x,y)\ll x^{-1}$, by \eqref{eq:1}, we see that Lemma \ref{prop:1} is
trivial when $x> 1$. We therefore assume that $x\leq 1$ for the
remainder of this section. 

Since $f\in \cW_1(V)$ there exists an absolute  constant $0<A\ll 1$ such that 
$f(v)=0$ unless
$v=(v^{(1)},\ldots,v^{(r_1+r_2)})$ satisfies $|v^{(l)}|\leq A$ for $1\leq l\leq r_1+r_2$.
We introduce a parameter $T$, to be selected in due course, which satisfies $0<T\leq A$.
We may then write
$$
I(x)=\sum_{S\subseteq \{1,\ldots,r_1+r_2\}} \int_{V_S} f(v) h(x,\nm(v))\d v,
$$
where $S$ runs over all subsets of $\{1,\ldots,r_1+r_2\}$ and 
$V_S$ is the set of $v\in V$ for which 
$$
|v^{(l)}|\begin{cases}\leq T, & \mbox{if $l\in S$,}\\
> T, & \mbox{if $l\not\in S$.}
\end{cases}
$$
Let us denote by $I_S(x)$ the integral over $V_S$.
In our proof of Lemma \ref{prop:1} we will obtain an asymptotic formula for $I_{\{1,\ldots,r_1+r_2\}}(x)$ and upper bounds for all other $I_S(x)$.

Given 
a vector $\ma{e}\in \ZZ_{\geq 0}^{k}$ and a vector
$\ma{t}=(t_1,\ldots t_k)$, we write
$
\ma{t}^{\ma{e}}=t_1^{e_1}\cdots t_{k}^{e_{k}}
$
for the associated monomial of degree 
$$
|\ma{e}|=\sum_{1\leq i\leq k} e_i.
$$
Let us work with a particular set $S$.
By symmetry we may assume without loss of generality that 
$S=\{1,\ldots,m\}\cup \{r_1+1,\ldots, r_1+n\}$ for some $m\leq r_1$ and $n\leq r_2$. 
Throughout this section it will be convenient to work with the coordinates 
$v=(\v,\y+i\z)$ on $V$, and to put
$$
\v=(\v',\v''), \quad 
\y=(\y',\y''), \quad \z=(\z',\z''),
$$
where $\v'=(v_1,\ldots,v_m)$ and $\v''=(v_{m+1},\ldots,v_{r_1})$, and similarly for $\y,\z$.
The idea is now to take a power series expansion of $f$ around the point 
$(\0,\v'',\0,\y''+i\z'')$ 
to produce an approximating polynomial of degree $M$ in the variables $\v',\y',\z'$, with error $O_M(\lambda_f^{M+1}T^{M+1})$. This leads to an expression
\begin{equation}\label{eq:taylor-f}
f(v)=
\sum_{\substack{
\ma{a} \in \ZZ_{\geq 0}^{m}, \ma{b},\ma{c}
\in \ZZ_{\geq 0}^{n}
\\
|\ma{a}|+|\ma{b}|+|\ma{c}|\leq M
}}
\phi_{\ma{a},\ma{b},\ma{c}}
\v'^{\a}\y'^{\b}\z'^{\c}
+O_M\left( 
\lambda_f^{M+1}T^{M+1}\right),
\end{equation}
for suitable coefficients
$\phi_{\ma{a},\ma{b},\ma{c}}=
\phi_{\ma{a},\ma{b},\ma{c}}(\v'',\y'',\z'')$,
with $|\phi_{\ma{a},\ma{b},\ma{c}}|=O_M(\lambda_f^M)$
and 
$$
\phi_{\ma{0},\ma{0},\ma{0}}=f(\ma{0},\v'',\ma{0},\y''+i\z'').
$$
Note that \eqref{eq:taylor-f} holds trivially if $S=\emptyset$, since then $f(v)=\phi_{\ma{0},\ma{0},\ma{0}}$.

Next we set 
$
B(H)=\{(y,z)\in \RR^2: y^2+z^2\leq H^2\}
$
for the ball of radius $H$ centred on the origin.
We may therefore write
\begin{equation}\label{eq:roo}
I_S(x)=
\sum_{\substack{
\ma{a} \in \ZZ_{\geq 0}^{m}, \ma{b},\ma{c}
\in \ZZ_{\geq 0}^{n}
\\
|\ma{a}|+|\ma{b}|+|\ma{c}|\leq M
}}
\int_{\v'',\y'',\z''} 
\phi_{\ma{a},\ma{b},\ma{c}}  
K
\d \v''\d \y'' \d\z''
+O_M\left( 
\frac{\lambda_f^{M+1}T^{M+1}}{x}\right),
\end{equation}
the integral being  over $\v''\in ([-A,T]\cup[T,A])^{r_1-m}$ and 
$(\y'',\z'')\in (B(A)\setminus B(T))^{r_2-n}$, and where
$$
K=
\int_{[-T,T]^m}\int_{B(T)^n}
\v'^{\a}\y'^{\b}\z'^{\c}
h(x,\nm(\v,\y+i\z)) \d \v'\d \y' \d\z'.
$$
Recalling the definition \eqref{eq:norm-trace} of $\nm(v)$,  we
see that
$$
\nm (\v,\y+i\z)=v_{1}\cdots v_{r_1}(y_1^2+z_1^2)\cdots (y_{r_2}^2+z_{r_2}^2).
$$
We extend this to $\nm (\v',\y'+i\z')$ and $\Phi=\nm (\v'',\y''+i\z'')$ in the obvious way, so that 
$\nm (\v,\y+i\z)=\Phi\nm (\v',\y'+i\z')$. We now make the change of variables 
$$
(\tilde \v,\tilde \y,\tilde \z)= |\Phi|^{1/(m+2n)} (\v',\y',\z'),
$$
which leads to the conclusion
\begin{equation}\label{eq:K}
K=|\Phi|^{-(|\a|+|\b|+|\c|)/(m+2n)-1} 
I_{m,n}^{(\a,\b,\c)} (x; T|\Phi|^{1/(m+2n)}),
\end{equation}
where for any $Y>0$ we set
\begin{equation}
\label{eq:koala}
I_{m,n}^{(\a,\b,\c)} (x; Y)
=
\int_{[-Y,Y]^m}\int_{B(Y)^n}
{\tilde \v}^{\a}\tilde{\y}^{\b}\tilde{\z}^{\c}
h(x,\nm(\tilde \v,\tilde \y+i\tilde\z)) \d \tilde\v \d \tilde\y \d\tilde\z.
\end{equation}
Note that $I_{m,n}^{(\a,\b,\c)} (x; Y)
$ is completely independent of the variables $\v'',\y'',\z''$. The analysis of this integral is rather involved and will lead to the following estimate.

\begin{lemma}\label{lem:I1-ij}
Let $N_1,N_2>0$, let $0<Y\ll 1$ and let $m,n\in \ZZ_{\geq 0}$, with $m\leq r_1$ and $n\leq r_2$. Let 
$\ma{a} \in \ZZ_{\geq 0}^{m}$ and $\ma{b},\ma{c}
\in \ZZ_{\geq 0}^{n}$, with 
$|\ma{a}|+|\ma{b}|+|\ma{c}|\leq M$.
Then we have 
$
I_{m,n}^{(\ma{a},\ma{b},\ma{c})}(x;Y) =0$ unless $a_i,b_j,c_j$ are all even,
in which case
$$
I_{m,n}^{(\ma{a},\ma{b},\ma{c})}(x;Y) =
-\frac{2^{r_1}\pi^{r_2}
  \Upsilon_K c_{m,n}^{(\ma{a},\ma{b},\ma{c})} }{\Delta_K} 
+
O_{M,N_1} \left( \left\{\frac{x}{Y^{m+2n}}\right\}^{N_1}\right)
+
O_{N_2}(x^{N_2}),
$$
with $\Upsilon_K$ given by \eqref{eq:Psi} and 
$$
c_{m,n}^{(\ma{a},\ma{b},\ma{c})}
=\begin{cases}
1, & \mbox{if $(m,n)=(r_1,r_2)$ and $
(\ma{a},\ma{b},\ma{c})=(\ma{0},\ma{0},\ma{0}),$} \\
0, & \mbox{otherwise}.
\end{cases}
$$
\end{lemma}

Let us delay the proof of Lemma \ref{lem:I1-ij} momentarily, in order to see how it can be used to conclude the proof of Lemma \ref{prop:1}. 
It is clear that $|\Phi|\geq T^{d-(m+2n)}$ for any $\v'',\y'',\z''$ appearing in \eqref{eq:roo}. 
It therefore follows from inserting Lemma \ref{lem:I1-ij} into \eqref{eq:K} that 
\begin{align*}
K
&=|\Phi|^{-(|\a|+|\b|+|\c|)/(m+2n)-1} 
\left\{
-\frac{2^{r_1}\pi^{r_2}
  \Upsilon_K c_{m,n}^{(\ma{a},\ma{b},\ma{c})}}{\Delta_K}
+
O_{M,N_1} \left( \left\{\frac{x}{T^d}\right\}^{N_1}\right)
+
O_{N_2}(x^{N_2})\right\}.
\end{align*}
Here we have 
$|\Phi|^{-(|\a|+|\b|+|\c|)/(m+2n)-1} \leq T^{-d(M/(m+2n)+1)}\leq T^{-d(M+1)}$.
Noting that $\Phi=1$ when $(m,n)=(r_1,r_2)$, 
and integrating trivially over $\v'',\y'',\z''$, 
we deduce from 
 \eqref{eq:roo} that
\begin{align*}
I_S(x)=~&
-\frac{2^{r_1}\pi^{r_2}
  \Upsilon_K c_{m,n}^{(\ma{0},\ma{0},\ma{0})}}{\Delta_K} 
f(0)
+O_M\left( 
\frac{\lambda_f^{M+1}T^{M+1}}{x}\right)\\
&+
O_{M,N_1,N_2}\left(
\lambda_f^M 
T^{-d(M+1)}\left(
\left\{\frac{x}{T^{d}}\right\}^{N_1}
+x^{N_2} 
\right)
\right).
\end{align*}
Let $N>0$ be an integer.
We now make the selection 
$
T=x^{1/(2d)},
$
which is clearly $O(1)$.
We choose $M$ so that $M+1=2d(N+1)$. This ensures that the first error term is satisfactory for Lemma \ref{prop:1}.
The second error term is seen to be satisfactory on choosing $N_1$ and $N_2$ sufficiently large in terms of $M$.
It remains to deduce from \eqref{eq:Delta} and \eqref{eq:Psi} that 
$$
-\frac{2^{r_1}\pi^{r_2}
 \Upsilon_K}{\Delta_K}= \frac{\sqrt{D_K}}{2^{r_2}}.
$$
Taken together,  we may now conclude that 
\begin{align*}
I(x)
&=I_{\{1,\ldots,r_1+r_2\}}(x)+\sum_{S\subsetneq \{1,\ldots,r_1+r_2\}} I_S(x)\\
&=
\frac{\sqrt{D_K}}{2^{r_2}} \cdot f(0) +O_N(\lambda_f^{2d(N+1)} x^N),
\end{align*}
as required to complete the proof of Lemma \ref{prop:1}.

\begin{proof}[Proof of Lemma \ref{lem:I1-ij}]
For notational convenience let us write $(\tilde \v,\tilde \y,\tilde \z)=(\v,\y,\z)$ in 
the integral $I_{m,n}^{(\a,\b,\c)} (x; Y)$ defined in 
\eqref{eq:koala}.
To begin with,  \eqref{eq:Hxy} yields
\begin{align*}
I_{m,n}^{(\ma{a},\ma{b},\ma{c})}(x;Y) 
&=x^{-1}
\int_{[-Y,Y]^{m}} \int_{B(Y)^{n}}  \v^{\ma{a}} \y^{\ma{b}}\z^{\ma{c}}
H\left(\frac{
|\nm(\v,\y+i\z)|
}{x}\right)
 \d \v \d \y \d \z
+O_{N_2}(x^{N_2}),
\end{align*}
for any $N_2>0$.
Note that the main term vanishes unless all the components of 
$\ma{a}, \ma{b}$ and $\ma{c}$ are even,
which we now assume. Now it is clear that 
$$
\int_{[-Y,Y]^{m}} \v^{\ma{a}} \d \v= \frac{2^{m}Y^{m+|\ma{a}|}}{\prod_{1\leq i \leq m}
(a_i+1)}.
$$
Likewise we have 
$$
\int_{B(Y)^{n}} \y^{\ma{b}}\z^{\ma{c}}\d \y \d \z=
\prod_{1\leq j\leq n} 
F(b_j,c_j),
$$
with 
$$
F(b,c)=
\int_{B(Y)} y^{b}z^{c} \d y \d z,
$$
for  even integers $b$ and $c$.
Switching  to polar coordinates, we 
call upon the identities
found in Gradshteyn and Ryzhik \cite[\S 2.511]{g-r}, 
deducing that
$$
\int_0^{2\pi} \cos^{2p}(t) \sin^{2q}(t)  \d t = 
\frac{(2p-1)!! (2q-1)!!}{(2p+2q)!! } \cdot 2\pi,
$$
for any integers $p,q\geq 0$.
Hence it follows that
\begin{align*}
F(b,c)
&=\int_{0}^{Y} \int_0^{2\pi}
\rho^{b+c+1} 
 \cos^{b}(\theta) \sin^{c}(\theta)  \d \rho \d \theta\\
&=
\frac{Y^{b+c+2}}{b+c+2} \cdot 
\frac{(b-1)!! (c-1)!!}{(b+c)!! } \cdot 2\pi,
\end{align*}
whence 
$$
\int_{B(Y)^{n}} \y^{\ma{b}}\z^{\ma{c}}\d \y \d \z=
\frac{C_{\b,\c}Y^{2n+|\b|+|\c|}}{\prod_{1\leq j\leq n} (b_j+c_j+2)}
$$
with 
$$C_{\b,\c}=
(2\pi)^{n}
\prod_{1\leq j\leq n} 
\left(
\frac{(b_j-1)!! (c_j-1)!!}{(b_j+c_j)!! } \right).
$$

It will be convenient to define
\begin{equation}\label{eq:lambda}
\nu_{\ma{a},\ma{b},\ma{c}}=
2^{m}C_{\b,\c}
\prod_{1\leq i\leq m}
\left(
\frac{1}{a_i+1}
\right)
\prod_{1\leq j\leq n}
\left(
\frac{1}{b_j+c_j+2}
\right)
\int_{0}^{\infty} \frac{w(u)}{u}\d u.
\end{equation}
This allows us to conclude that
\begin{equation}\label{eq:friday}
I_{m,n}^{(\a,\b,\c)}(x;Y)=
\frac{\nu_{\ma{a},\ma{b},{\ma{c}}} Y^{|\a|+|\b|+|\c|+m+2n}}{x}
- 
\frac{2^{m}C_{\b,\c}}{\Delta_K}
\sum_{\fz} g(\n \fz)+  O_{N_2}(x^{N_2}),
\end{equation}
with 
\begin{align*}
g(t)
&=\frac{1}{C_{\b,\c}}
\int_{[0,Y]^{m}} \int_{B(Y)^{n}}  \v^{\ma{a}} \y^{\ma{b}}\z^{\ma{c}}
\cdot
\frac{1}{xt} w\left(
\frac{
|\nm(\v,\y+i\z)|
}{xt}\right)
 \d \v \d \y \d \z,
\end{align*}
for $t>0$.
Switching to polar coordinates as above we may write
\begin{align*}
g(t)
&=
\int_{[0,Y]^{m}} \int_{[0,Y]^{n}}  \v^{\ma{a}} 
\brh^{\b+\c+\mathbf{1}}
\cdot
\frac{1}{xt} w\left(
\frac{
v_1\cdots v_{m}\rho_1^2\cdots \rho_{n}^2
}{xt}\right)
 \d \v \d \brh.
\end{align*}
We are therefore led to analyse the sum 
$$
S=\sum_{\fz} g(\n \fz),
$$
for which we will use properties of the Mellin transform
$$
\hat g(s)=\int_0^\infty t^{s-1}g(t)\d t.
$$
It follows that 
\begin{align*}
S=\sum_{n=1}^\infty a_n g(n)
=\frac{1}{2\pi i} \int_{(2)} \zeta_K(s) \hat g(s) \d s,
\end{align*}
where $a_n$ are the coefficients appearing in the Dedekind zeta function $\zeta_K(s)$
and the integral is over the line $\sigma=\Re(s)=2$.
In order to  move the line of integration further to the left we will need a better understanding of the
analyticity of $\hat g$.

We  will write 
$
u=(xt)^{-1}v_1\cdots v_{m}\rho_1^2\cdots \rho_{n}^2
$
and substitute for one of the variables. 
If $m\neq 0$ then we substitute for $v_1$.  
Alternatively, if $m=0$, then we substitute for $\rho_1$ and argue similarly. 
Assuming without loss of generality that $m\neq 0$, 
it will be convenient to put 
$L=v_2\cdots v_{m}\rho_1^2\cdots \rho_{n}^2$ and to set
$\v_1=(v_2,\ldots,v_{m})$ and 
$\ma{a}_1=(a_2,\ldots,a_{m})$.
Then for  $\Re(s)>0$ it follows that
\begin{align*}
\hat g(s)
&=
\int_0^\infty t^{s-1}
\int_{[0,Y]^{m-1}} \int_{[0,Y]^{n}}  
\frac{\v_1^{\ma{a}_1} 
\brh^{\b+\c+\mathbf{1}}}{L}
\int_0^{\frac{LY}{xt}}
\left(\frac{uxt}{L}\right)^{a_1} w(u) \d u \d \v_1 \d\brh  \d t\\
&=
x^{a_1}
\int_0^\infty 
u^{a_1} w(u)
\int_{[0,Y]^{m-1}} \int_{[0,Y]^{n}}  
\frac{\v_1^{\ma{a}_1} \brh^{\ma{b}+\c+\mathbf{1}}}{L^{a_1+1}}
\int_0^{\frac{LY}{xu}}
t^{a_1+s-1}
 \d t \d \v_1 \d\brh  \d u.
\end{align*}
Carrying out the integration over $t$, 
we obtain
\begin{align*}
\hat g(s)
&=
\frac{Y^{a_1+s}}{(a_1+s)x^s}
\int_0^\infty 
\frac{ w(u)}{u^{s}} \d u
\int_{[0,Y]^{m-1}} \int_{[0,Y]^{n}}  
\v_1^{\ma{a}_1} \brh^{\ma{b}+\ma{c}+\mathbf{1}} L^{s-1}
 \d \v_1 \d\brh= F_Y(s),
\end{align*}
with
$$
F_R(s)
=
\frac{ 
R^{|\a|+|\b|+|\c|+(m+2n)s}}{x^s
\prod_{1\leq i\leq m}(a_i+s)
\prod_{1\leq j\leq n}(b_j+c_j+2s)}
\int_0^\infty 
\frac{ w(u)}{u^{s}} \d u,
$$
for $\Re(s)>0$ and $R>0$. 
As a formula this continues to make sense in the
half-plane $\Re(s)\leq 0$ and so provides a meromorphic continuation of
$\hat g$ to the whole of $\CC$, with at most  poles
at the non-positive integers.
One sees that the pole at $s=0$ has order at most $m+n-1$ when 
$(\ma{a},\ma{b},\ma{c})\neq (\ma{0},\ma{0},\ma{0})$ and  order 
$m+n$ otherwise. Likewise, the poles at the negative even (resp.\ odd) integers have  order at most $m+n$ (resp.\ $n$).

We take this opportunity to record an upper
bound for $|F_R(s)|$. Assume that $s=\sigma+it$ with $|t|\geq 1$ and recall that $w\in
\cW_1^+(\RR)$. 
Repeated integration by parts then yields
$$
\int_0^\infty 
\frac{ w(u)}{u^{s}} \d u=
\frac{1}{(s-1)\cdots (s-N)}
\int_0^\infty 
u^{N-s} w^{(N)}(u) \d u,
$$
for any integer $N\geq 0$. In this way we conclude that 
$$
|F_R(s)|
 \ll_{M,N,\sigma}  \frac{R^{|\a|+|\b|+|\c|}}{(1+|t|)^N}
\left(\frac{R^{m+2n}}{x}\right)^{\sigma},
$$
if $\Im (s)\geq 1$,  for any $N\geq 0$.

Returning to our formula for $S$ we seek to move the line of
integration back to the line $\sigma=-N_1$ for some positive constant
$N_1.$ This is facilitated by the polynomial decay in $|t|$  that we have observed in $F_R(s)$.
Indeed, in view of our convexity estimates \eqref{eq:mu}, we are able to shift the line of
integration arbitrarily far to the left. 
In doing so we will encounter poles at $s=1$ and possibly also at the non-positive  integers. 
We note that if $\hat g(s)$ has a pole of order $\leq r_1+r_2-1$ at $s=0$ 
then it will be compensated for by the presence of $\zeta_K(s)$, which has a zero 
of order $r_1+r_2-1$ at $s=0$. 
Similarly, any poles at the negative even (resp.\ odd) integers will be compensated for by the zeros of $\zeta_K(s)$ 
of order $r_1+r_2$ (resp.\ $r_2$)
at these places.
The residue at $s=1$ of $\zeta_K(s)\hat g(s)$ is 
$$
 \Delta_K\hat g(1) = \frac{\Delta_K 
}{2^{m} C_{\b,\c}} 
\cdot 
\frac{ \nu_{\a,\b,\c}
Y^{|\a|+|\b|+|\c|+m+2n}}{x},
$$
in the notation of \eqref{eq:Delta} and \eqref{eq:lambda}.
When 
$(\ma{a},\ma{b},\ma{c})=(\ma{0},\ma{0},\ma{0})$ and $(m,n)=(r_1,r_2)$,
then $\zeta_K(s) \hat g(s)$ has  a simple pole at $s=0$ with residue
$$
\frac{\Upsilon_K}{2^{r_2}}  \int_{0}^{\infty} w(u)\d u = 
\frac{\Upsilon_K}{2^{r_2}},
$$
where
$\Upsilon_K$ is given by \eqref{eq:Psi}.  Putting this together we may therefore conclude that 
\begin{align*}
S=~&
\frac{\Delta_K 
}{2^{m}C_{\b,\c}} 
\cdot 
\frac{
 \nu_{\a,\b,\c}
Y^{|\a|+|\b|+|\c|+m+2n}}{x}
+
\frac{\Upsilon_K c^{(\a,\b,\c)}_{m,n}}{2^{r_2}}\\
&+
O_{M,N_1} \left( Y^{|\a|+|\b|+|\c|}\left\{\frac{Y^{m+2n}}{x}\right\}^{-N_1}
\right),
\end{align*}
with $c^{(\a,\b,\c)}_{m,n}$ as in the statement of the lemma. 
Recalling that $Y\ll 1$ and 
substituting this into \eqref{eq:friday}, 
this therefore concludes the proof of Lemma \ref{lem:I1-ij}.
\end{proof}

\section{Application to hypersurfaces}\label{s:poisson-forms}

Suppose that $Y\subseteq \AA_K^{n}$ is a hypersurface defined by  
a  polynomial  $F\in \fo[X_{1},\ldots,X_{n}]$. In order to  
gauge whether or not $Y(\fo)$ 
is empty it is sometimes fruitful to study the asymptotic behaviour of sums
$$
 N_{W}(FP) = 
 \sum_{\x\in \fo^n} \delta_{K}(F(\x)) W(\x/P) ,
$$
as $P\rightarrow \infty$, where $W\in \cW_n(V)$.
For any $Q\geq 1$ we deduce from Theorem \ref{t:delta} that 
\begin{align*}
\No 
&=\frac{c_{Q}}{Q^{2d}}\sum_\mathfrak{b} \sigstar
\sum_{\x\in \fo^{n}}
   \sigma(F(\x))W(\x/P)h\left(\frac{\n\mathfrak{b}}{Q^d}  ,\frac{|\nm(F(\x))|}{Q^{2d}}\right).
\end{align*}
In view of the  fact that $h(x,y)\neq 0$ only if $x\leq \max(1,2|y|) $, 
it is clear that the sum over $\mathfrak{b}  $ is restricted to $\n \fb\ll Q^{d} $, if $Q$ is taken to be of order $P^{(\deg F)/2}$.
Breaking the inner sum over $\x$ into residue classes modulo $\fb$, we see that 
it can be written
$$
\sum_{\ma{a}\in (\oh/\mathfrak{b}  )^n}\sigma(F(\ma{a}))
\sum_{\x\in \mathfrak{b}^n}W\left((\x+\ma{a})/P\right)
h\left(\frac{\n \mathfrak{b}}{Q^d}  ,\frac{|\nm(F(\x+\ma{a}))|}{Q^{2d}}\right),
$$
for any primitive character $\sigma$ modulo $\fb$. 
We apply the usual multi-dimensional Poisson summation formula (in the form 
\cite[\S 5]{skinner}, for example), finding that 
the inner sum over $\x$ is 
\begin{align*}
\frac{2^{r_2 n}}{D_K^{n/2}(\n \mathfrak{b}  )^{n}} 
\sum_{\m\in \hat{\mathfrak{b}  }^n}
\e(\m.\a)
\int_{V^{n}} W(\x/P)
h\left(\frac{\n \mathfrak{b}}{Q^d}, \frac{|\nm(F(\x))|}{Q^{2d}}\right)
\e(-\m.\x)\d \x,
\end{align*}
where $D_K$ is the absolute discriminant of $K$ and $\hat \fb$ is the dual of $\fb$ taken with respect to the trace.  
Putting everything together we have therefore established the following result.

\begin{theorem}
 \label{prop:2}
We have
$$
\No=
\frac{c_{Q}2^{r_2 n}}{D_{K}^{n/2}Q^{2d}}
\sum_{\mathfrak{b}} \sum_{\m\in \hat{\mathfrak{b}  }^n}(\n \mathfrak{b}  )^{-n}S_\mathfrak{b}  (\m)
I_\mathfrak{b}  (\m),
$$
where the sum over $\fb$ is over non-zero integral ideals and
\begin{align*}
S_\mathfrak{b}  (\m)
&=\sigstar\sum_{\ma{a}\bmod{\mathfrak{b} }}\sigma(F(\ma{a}))\e(\m. \ma{a}),\\
I_\mathfrak{b}  (\m)
&=\int_{ V^{n}}W(\x/P)h\left(\frac{\n \mathfrak{b}}{Q^d}  ,\frac{\nm(F(\x))}{Q^{2d}}\right)\e\left(-\m.\x \right)\d\x.
\end{align*}
\end{theorem}

This result is a number field analogue of \cite[Thm.~2]{H}.
We apply this in the case that 
$F\in \fo[X_1,\ldots,X_n]$ is 
a non-singular cubic form in $n\geq 10$ variables.
Using the embedding of $K$ into $V$, we may write 
$$
F=\bigoplus_{l=1}^{r_1+r_2} F^{(l)},
$$ 
where each $F^{(l)}$ is a
cubic form over  $K_l$ for $1\leq l\leq r_1+r_2$.
Recall the definition of the norms $\langle\cdot \rangle, \|\cdot\|$ from \S \ref{s:technical}.

Let 
$
\bxi=\bigoplus_l \bxi^{(l)} \in V^n
$ 
be a suitable point chosen as in \cite[Lemma 13(i)]{skinner}. 
Let $\delta_0>0$ be a small 
constant such that the
inverse function theorem can be used for $F$ in the region $\langle \x-\bxi\rangle\leq \delta_0$.
Let $\u=\x-\bxi$. Thus,  without loss of generality, we may assume that there exists a smooth function $f$, 
such that for each $1\leq l\leq r_1+r_2 $ we have 
\begin{equation}\label{eq:space}
u^{(l)}_1=f^{(l)}(F^{(l)}(\u^{(l)}+\bxi^{(l)}),u_2^{(l)},\ldots,u_n^{(l)}),
\end{equation}
whenever  $|\u^{(l)}|\leq \langle \u\rangle \leq \delta_0$.
We henceforth view $\delta_0$ as being fixed once and for all.

Next, we take 
$w_0\in \cW_n^+(V)$ to be a smooth weight function which takes the value $1$ on the region $\langle \x-\bxi\rangle\leq
\delta_0/2$ and is zero
outside the region $\langle \x-\bxi\rangle\leq \delta_0$.
Let $\omega:V^n\rightarrow \RR_{\geq 0}$ be the smooth weight function
\begin{equation}\label{eq:omega}
\omega(\x)=\exp(-(\log P)^4\|\x-\bxi \|^2).
\end{equation}
Note that $\omega(\x)$ is very small unless $\langle \x-\bxi\rangle\ll 1/(\log P)^2$.
The function $W:V^n\rightarrow \RR_{\geq 0}$ that we shall work with 
in   $N_{W}(F,P)$
is
$$
W(\x)=w_0(\x)\omega(\x).
$$
In particular $W$ is supported on the region $\langle \x\rangle\ll 1$.
Let $Q=P^{3/2}$.
It now follows from Theorem \ref{prop:2} that
$$
\No=
\frac{c_{Q}2^{r_2 n}}{D_{K}^{n/2}Q^{2d}}
\sum_{\substack{(0)\neq \mathfrak{b}\subseteq \fo\\
\n \mathfrak{b}  \ll Q^{d}
}} \sum_{\m\in \hat{\mathfrak{b}  }^n}(\n \mathfrak{b}  )^{-n}S_\mathfrak{b}  (\m)
I_\mathfrak{b}  (\m).
$$
The definition of $h$ in 
$I_\mathfrak{b}  (\m)$
means that one can freely replace $|\nm(F(\x))|$ by $\nm(F(\x))$.

We will use this expression to obtain an asymptotic lower bound for $\No$, as $P\rightarrow \infty$.   The terms corresponding to $\m=\0$ will give the main contribution.
The following section will be concerned with the exponential integrals 
$I_\mathfrak{b}  (\m)$. 
Our analysis of the complete exponential sums $S_{\mathfrak{b}}(\m)$ will take place in \S \ref{s:sums}. Finally, in \S \ref{s:conclusion} we shall handle the terms with $\m=\0$ and  draw together the various estimates in order to conclude the proof of Theorem \ref{t:10}.

Let $\rho=Q^{-d}\n\bbb\ll 1$.
The first task in \S \ref{s:integrals} will be to prove the following result,
which is based on repeated integration by parts.

\begin{lemma}
\label{ibeta}
For any non-zero $\m\in  V^n$ and any integer $N\geq 0$, we have
$$
I_\fb( \m)\ll_{N} \rho^{-1}P^{dn} \left(\frac{(\log P)^2}{
\rho P \langle \m\rangle }\right)^{N}. 
$$
\end{lemma}

This result shows that 
$I_{\fb}(\m) $  decays faster than any polynomial decay
in $\langle \m\rangle$.   
Noting the trivial bound 
$|S_{\fb}(\m)|\leq (\n \fb)^{n+1}$,
the tail of the series involving $\m$ in our expression for 
$\No$  therefore makes 
a negligible contribution,   leaving us free to truncate the sum over $\m$ by $\langle \m\rangle \leq  P^A$,  for some appropriate absolute constant $A>0$.
Thus we have 
\begin{equation} \label{eq:no}
\No=
\frac{c_{Q}2^{r_2 n}}{D_{K}^{n/2}Q^{2d}}
\sum_{\substack{(0)\neq \mathfrak{b}\subseteq \fo\\
\n \mathfrak{b}  \ll Q^{d}
}} \sum_{
\substack{\m\in \hat{\mathfrak{b}  }^n\\
\langle \m \rangle \leq P^A}}
(\n \mathfrak{b}  )^{-n}S_\mathfrak{b}  (\m)
I_\mathfrak{b}  (\m) +O(1).
\end{equation}

We close this section by recording some further notation that will feature  in our analysis.
For a subset $S\subset \{1,\ldots,r_1+r_2\}$, we define the restricted norm
$$
\nm_S(v)=\prod\limits_{\substack{l\in S\\ l\leq
r_1}}v^{(l)}\prod\limits_{\substack{l\in S\\ l>r_1}}|v^{(l)}|^2,
$$ 
on $V$. We follow the convention that 
$\nm_{\emptyset}(v)=1$ and define the ``degree of $S$'' to  mean
$d(S)=\sum_{l\in S}c_l$, with 
$c_l$ given by \eqref{eq:cl}.
For any $v\in V$, let 
\begin{equation}\label{eq:Tv}
T(v)=\{1\leq l\leq r_1+r_2: |v^{(l)}|>1\}.
\end{equation}
The quantity $|\nm_{T(v)}(v)|$ will provide a useful measure of the ``height'' of a point $v\in V$, and we henceforth set
\begin{equation}\label{eq:Tv'}
\cH: V\rightarrow \RR_{>0}, \quad 
v\mapsto |\nm_{T(v)}(v)|.
\end{equation}
We may now establish the following result.

\begin{lemma}
Let $A>0$ and let $\alpha<-1$. Then we have 
$$
\int_{\substack{\{v\in V:
\cH(v)\leq A\}}} 
\cH(v)^{\alpha}\d v\ll_\alpha 1,
$$
uniformly in $A$.
\label{nmintegral}
\end{lemma}

\begin{proof}
Given any subset $S$ of $\{1,\ldots,r_1+r_2\}$, let $R_S$ denote the set of $v\in V$ for which $|\nm_S(v)|\leq A$ and $T(v)=S$. In order to establish the  lemma it suffices to
prove the desired bound for each integral 
$$
I_S=\int_{R_S} \cH(v)^{\alpha}\d v.
$$
Note that for $T(v)=\emptyset$, we have $\cH(v)=1$, whence 
$I_\emptyset\ll 1$ in this case. We next assume  that 
$S=\{s_1,...,s_{l+m}\}$ is non-empty, with 
$s_i\leq r_1$ for $1\leq i\leq l$ and 
$s_{l+i}>r_1$ for $1\leq i\leq m$.
Let us make the polar change of variables 
$(v^{(s_1)}, \ldots, v^{(s_{l+m})})$ goes to 
$(u_1,\ldots,u_{l+m},\theta_1, \ldots, \theta_{m})$, with 
$$u_i=
\begin{cases}
		|v^{(s_i)}|,   & \mbox{if $1\leq i\leq l$,} \\
		|v^{(s_{i})}|^2, & \mbox{if $l<i\leq l+m$.}
	\end{cases}
$$
In particular 
$\d v^{(s_{l+i})}=2^{-1}\d u_{l+i}\d \theta_i$, for $1\leq i\leq m$ and  
$\cH(v)=u_1\cdots u_{l+m}$.
It follows that 
$$
I_S\ll \int_{\substack{
u_1,\ldots,u_{l+m}\geq 1\\
u_1\cdots u_{l+m}\leq A}}  (u_1\cdots u_{l+m})^{\alpha} \d u_1\cdots \d u_{l+m}.
$$
The statement of the lemma is clearly trivial unless $A\geq 1$, which we now assume.
Write $\sigma=l+m$ for the cardinality of $S$ and 
denote by $J_{\sigma}$ the integral on the right hand side.
In order to complete the proof of the lemma 
it suffices to show that 
$J_\sigma\ll_{\alpha} 1$.
We do so by induction on $\sigma$, the case $\sigma=1$ being trivial.
For $\sigma>1$ we integrate over $u_\sigma$, finding that 
$
J_\sigma \ll_\alpha J_{\sigma-1}\ll_\alpha 1,
$
by the induction hypothesis.
This completes the proof of the lemma.
\end{proof}

Next, for any $v\in V$ let
\begin{equation}\label{eq:def-I}
\mathfrak{I}(v)=
\int_{\substack{\x\in V^n\\ F(\x)=v}} W(\x)\d\x,
\end{equation}
with $W$ as above.
It is easy to see that 
$\mathfrak{I}(v)$ 
is compactly supported. 
We claim that it is also an infinitely differentiable function on $V$.
To see this we note first that
\begin{equation}\label{eq:black}
\begin{split}
 \wtil(v)
&=\prod_{l=1}^{r_1+r_2}\int_{
\substack{\x^{(l)}\in K_l^n\\ 
F^{(l)}(\x^{(l)})=v^{(l)}}}
W(\x^{(l)})\d\x^{(l)}\\
&=\prod_{l=1}^{r_1+r_2}
 \wtil^{(l)}(v^{(l)}),
\end{split}
\end{equation}
say.
We need to  show that 
$ \wtil^{(l)}(v^{(l)})$ is infinitely differentiable function on $K_l$, for each choice of $1\leq l\leq r_1+r_2$. 
We will give details for the case $l\leq r_1$ only. The case $l>r_1$ follows in the same way (see \cite[page 464]{skinner} for a similar calculation).
Our first step is to make the change of variables 
$\u^{(l)}=\x^{(l)}-\bxi^{(l)}$,  writing
$\widetilde W(\u^{(l)})=W(\u^{(l)}+\bxi^{(l)})$ for convenience of notation. 
According to \eqref{eq:space} we have $u^{(l)}_1=f^{(l)}(v^{(l)},u_2^{(l)},\ldots, u_n^{(l)})$, 
whence 
$$
 \wtil^{(l)}(v^{(l)})=\int_{\RR^{n-1}} \partial_1 f^{(l)}(v^{(l)},\z^{(l)}) \widetilde W\left(f^{(l)}(v^{(l)},\z^{(l)}),\z^{(l)} \right) \d\z^{(l)},
$$
where $\z^{(l)}=(u_2^{(l)},\ldots, u_n^{(l)})$ and $\partial_1 f^{(l)}$ denotes the derivative with respect to the first coordinate.
Since $f^{(l)}$ is smooth
this shows that 
$ \wtil^{(l)}(v^{(l)})$ is infinitely differentiable on $\RR$. In fact, 
for any $N\in \ZZ_{\geq 0}$, 
the $N$th derivative of 
$\partial_1 f^{(l)}(v^{(l)},\z^{(l)})$ with respect to $v^{(l)}$ 
is $O_N(1)$ on the support of $\widetilde W$.
Moreover,  it is easy to see that the derivatives of $\omega(\x^{(l)})$ are bounded by $O_N((\log P)^{2N}).$ 
In the notation of \eqref{eq:sobolev}, 
we may therefore conclude that 
\begin{equation}\label{eq:r1}
\lambda^N_{\mathfrak{I}} \ll_N (\log P)^{2N},
\end{equation}
for any $N\in \ZZ_{\geq 0}$.

It turns out that 
$\mathfrak{I}(0)$ is the ``singular integral'' for the problem. The ``singular series'' is formally given by the  infinite sum 
\begin{equation} \label{eq:sing-series}
\fS=
\sum_{(0)\neq \mathfrak{b}\subseteq \fo}
(\n \mathfrak{b}  )^{-n}S_\mathfrak{b}  (\0).
\end{equation}
Our main aim is to establish the following result, which clearly suffices for  Theorem \ref{t:10}.

\begin{theorem}\label{t:10'}
Assume that $n\geq 10$. Then there exists
$\Delta>0$ such that 
$$
\No=
\frac{2^{r_2 (n-1)}}{D_{K}^{(n-1)/2}}
\fS \mathfrak{I}(0)P^{(n-3)d}+O(P^{(n-3)d-\Delta}), 
$$
with 
$$
(\log P)^{-2d(n-1)} \ll
\fS\mathfrak{I}(0) \ll (\log P)^{-2d(n-1)}.
$$
\end{theorem}

\section{Cubic exponential integrals}\label{s:integrals}

Let $\fb$ be a non-zero integral ideal with  $\n\bbb\ll Q^d$ and define
$$
\rho=Q^{-d}\n\bbb.
$$
Thus  $\rho\in \RR$ satisfies $\rho\ll 1$.  In this section we shall produce a number of 
estimates for the exponential integral $I_{\fb}(\m)$ in \eqref{eq:no}, beginning with a proof of Lemma \ref{ibeta}.
By a change of variables we get
\begin{equation}\label{eq:exp-int'}
I_\fb( \m)=P^{dn}\int_{V^n} W(\x)h(\rho ,\nm(F(\x)))\e(-P\m.\x)\d\x,
\end{equation}
where $W$ is supported on the region $\langle \x\rangle\ll 1$. 
Note that
\begin{align*}
\partial^{\bbe} \left\{W(\x)h(\rho ,\nm(F(\x)))\right\}
&=\sum_{\bbe=\bbe_1+\bbe_2}\partial^{\bbe_1}W(\x)\partial^{\bbe_2}h(\rho,\nm(F(\x))),
\end{align*}
for $\bbe, \bbe_1, \bbe_2$ running over  $\mathbb Z_{\geq 0}^{dn}$.
Since $F$ is a  polynomial, there exist polynomials $f_{\bbe_2,j}$, such that 
\begin{equation*}
\partial^{\bbe_2} h(\rho,\nm(F(\x)))=\sum_{j=0}^{|\bbe_2|} f_{\bbe_2,j}(\x)h^{(j)}(\rho,\nm(F(\x))),
\end{equation*}
with $h^{(j)}(x,y)=\frac{\partial^j}{\partial y^j} h(x,y)\ll_j  x^{-j-1}$, by Lemma \ref{lem:decay}.
Moreover, 
$|\partial^{\bbe_1}
W(\x)|\ll_{|\bbe_1|} (\log
P)^{2|\bbe_1|}$  for any $\x\in V^n$.  
Combining these estimates,  together  with the inequality $\rho\ll 1$, we get
$$
\partial^{\bbe} \left\{ W(\x)h(\rho ,\nm(F(\x)))\right\}\ll_{|\bbe|} \rho^{-1}(\rho^{-1}(\log P)^2)^{|\bbe|}.
$$
Using repeated integration by parts in \eqref{eq:exp-int'}, we arrive at the statement of 
Lemma \ref{ibeta}.

\medskip

We now turn to a  more sophisticated treatment of $I_{\fb}(\m)$.
Let $L$ be a constant, with  $1\ll L\ll 1$,  such that $\langle F(\x)\rangle \leq L$ for every $\x \in 
\supp(W)$. 
Set
$w_2(v)=w_1(v/2L)$, where $w_1\in W_1^+(V)$ is any weight function  
which takes the value $1$ in the region  $\langle
v\rangle \leq 1$. 
Then we clearly have $W(\x)=W(\x)w_2(F(\x))$ for every $\x\in V^n$.
This allows us to write 
\begin{align*}
I_\fb( \m)
&= \int_{V^n} W(\x/P)\left\{w_2(Q^{-2}F(\x))h(\rho,\nm(Q^{-2}F(\x)))\right\}\e(-\m.\x)\d\x.
\end{align*} 
The  Fourier inversion formula implies that 
$$
w_2(Q^{-2}F(\x))h(\rho,\nm(Q^{-2}F(\x)))=
\int_V p_\rho (v) \e(vQ^{-2}F(\x))\d v, 
$$
where 
\begin{equation}\label{eq:p}
p_\rho (v)=\int_{V} w_2(x)h(\rho,\nm(x))\e(-vx)\d x. 
\end{equation}
This yields
\begin{equation}\label{ibc}
I_\fb(\m)
=\int_V p_\rho
(v)
K(Q^{-2}v,\m) \d v, 
\end{equation}
with 
$$
K(v,\m)=\int_{V^n} W(\x/P)\e(vF(\x)-\m.\x)\d\x.
$$

\subsection{Weighted exponential integrals}

It is clear from our work above that we will need good estimates for integrals of the form
$$
\int_{V} w(x)h(\rho,\nm(x))\e(-vx)\d x,
$$
for $w\in \cW_1^+(V)$.
In fact, in the context of \eqref{ibc}, in order for the integral over $v$ to converge 
we require an estimate for $p_\rho(v)$ that decays sufficiently fast.
It turns out that we will require savings of the form $\cH(v)^{-N}$, in the notation of \eqref{eq:Tv'}.  
The principal means of achieving this will be 
the use of   integration by parts repeatedly.
However, it will be crucial to apply this process in multiple directions, with respect to every  component  $v^{(l)}$ of $v$ such that $l\in T(v)$.

Recall our decomposition 
$v=(\v,\y+i\z)$, from \S \ref{s:weights}.
Let $\partial_l$ denote $\partial_{v_l}$ for $1\leq l\leq r_1$. For $l>r_1$, we let
$\partial_l$ denote
$\partial^2_{y_{l-r_1}}+\partial_{z_{l-r_1}}^2$. 
Given a subset $S$ of $\{1,\ldots ,r_1+r_2\}$, we let
$$
\partial_S=\prod_{l\in
S}\partial_l,
$$
and we let $S^c=\{1,\ldots,r_1+r_2\}\setminus S$.
Recall the notation 
$h^{(m)}(x,y)=\frac{\partial^m}{\partial y^m} h(x,y)$ for any $m\in \ZZ_{\geq 0}$.
We are now ready to establish the following result.

\begin{lemma}
\label{lemma1}
Let $S\subseteq \{1,\ldots,r_1+r_2\}$ and let
$w\in \cW_1(V)$.  Then given any  $k_1, k_2\in \ZZ_{\geq 0}$, there
exist weight functions 
$w_{S,k_1,k_2}, 
w_{S,k_1,k_2}^{(m)}\in \cW_1(V)$, such that 
\begin{align*}
\partial_{S}
\big\{w(v)&\nm(v)^{k_1}h^{(k_2)}(\rho,\nm(v))\big\}
\\
=~&
w_{S,k_1,k_2}(v)\nm(v)^{k_1}h^{(k_2)}(\rho,\nm(v))\\
& +
\sum_{m=0}^{d(S)}
w_{S,k_1,k_2}^{(m)}(v)
\nm_{S^c}(v)\nm(v)^{\max\{0,k_1+m-1\}}h^{(k_2+m)}
(\rho,\nm(v)).
\end{align*} 
\end{lemma}

\begin{proof}
The proof will be given using induction on the cardinality of the set $S$. We will first prove the lemma when $\#S\leq 1$. Let $2\leq M\leq r_1+r_2$. 
Then,  assuming the result 
 to
be true for all subsets of $\{1,\ldots,r_1+r_2\}$,  with  cardinality at most $M-1$,  
we will decompose a given subset $S\subseteq \{1,\ldots,r_1+r_2\}$ of cardinality $M$ as 
$S=S_1\cup \{j\}$, with  $\#S_1= M-1$.  We will then use the
fact that
$\partial_S=\partial_j\partial_{S_1}$ and use the induction hypothesis accordingly.
Notice that for $1\leq j\leq r_1$ we have  
$$
\partial_j \nm(v)=\nm_{\{j\}^c}(v),
$$
and for  $r_1<j\leq r_2$,
$$\frac{
\partial}{\partial {y_j}} \nm(v)=2y_{j}\nm_{\{j\}^c}(v),\quad
\frac{\partial}{\partial {z_j}} \nm(v)=2z_{j}\nm_{\{j\}^c}(v).
$$
Let $K(v)=w(v)\nm(v)^{k_1}h^{(k_2)}(\rho,\nm(v))$.

The case $S=\emptyset$ is trivial. Suppose next that $S=\{j\}$. By symmetry it suffices to  deal with cases
when $S=\{1\}$ and
$S=\{r_1+1\}$. Suppose first that $j=1$. Then
\begin{align*}
\partial_S
K(v)
=~&\partial_{1}(w(v)\nm(v)^{
k_1}h^{(k_2)}(\rho ,\nm(v))\\
=~&(\partial_1
w(v))\nm(v)^{k_1} h^{(k_2)}(\rho,\nm(v))\\
&+k_1w(v)\nm_{S^c}(v)\nm(v)^{k_1-1}h^{(k_2)}(\rho,\nm(v))\\
&+w(v) \nm_{S^c}(v)\nm(v)^{k_1}h^{(k_2+1)}(\rho,\nm(v)).
\end{align*}
This is clearly of the required form, with  $w_{S,k_1,k_2}=\partial_1 w$.
Suppose next that  $S=\{r_1+1\}$, so that  $\partial_S=\partial_{y_1}^2+\partial_{z_1}^2$. Notice that
\begin{align*}
\partial_{y_1}^2 K(v)
=~& 
\partial_{y_1}\left(\partial_{y_1}\left(w(v)\nm(v)^{k_1}h^{(k_2)}(\rho,
\nm(v))\right)\right)\\
=~&
\partial_{y_1}\left((\partial_{y_1}
w(v))\nm(v)^{k_1}
h^{(k_2)}(\rho,\nm(v))\right)\\
&+2k_1\partial_{y_1}\left(y_1\nm_{S^c}(v)
w(v)\nm(v)^{k_1-1
}
h^{(k_2)}(\rho,
\nm(v))\right)\\
&+\partial_{y_1}\left(2y_1\nm_{S^c}(v)w(v)
\nm(v)^{k_1}h^{(k_2+1)}(\rho,\nm(v))\right).
\end{align*}
We  deal here only with the third term in this expression, the remaining two terms being of a similar ilk. Let $f(v)=2\nm_{S^c}(v)w(v)
\nm(v)^{k_1}h^{(k_2+1)}(\rho,\nm(v))$. We find that 
\begin{align*}
\partial_{y_1}\left(y_1f(v)\right)
=~&
2\nm_{S^c}(v)\left\{ w(v)+y_1\partial_{y_1}w(v)\right\}\nm(v)^{k_1}h^{(k_2+1)}(\rho,\nm(v))\\
&+4k_1y_1^2(\nm_{S^c}(v))^2w(v)\nm(v)^{k_1-1}h^{(k_2+1)}(\rho,\nm(v))\\
&+4y_1^2(\nm_{S^c}(v))^2w(v)\nm(v)^{k_1}h^{(k_2+2)}(\rho,\nm(v)).
\end{align*}
We carry out the same process for $z_1$ and get similar expressions. Adding together the expressions corresponding to
$y_1$ and $z_1$ we get the contribution
\begin{align*}
\partial_{y_1}\left(y_1f(v)\right)
&+
\partial_{z_1}\left(z_1f(v)\right)\\
=~&
2\nm_{S^c}(v)\left\{ 2w(v)+y_1\partial_{y_1}w(v)+z_1\partial_{z_1}w(v)\right\}
\nm(v)^{k_1}h^{(k_2+1)}(\rho,\nm(v))\\
&+4k_1(y_1^2+z_1^2)(\nm_{S^c}(v))^2w(v)
\nm(v)^{k_1-1}h^{(k_2+1)}(\rho,\nm(v))\\
&+4(y_1^2+z_1^2)(\nm_{S^c}(v))^2w(v)\nm(v)^{k_1}h^{(k_2+2)}(\rho,\nm(v))\\
=~&
E_1(v)\nm_{S^c}(v)\nm(v)^{k_1}h^{(k_2+1)}(\rho,\nm(v))\\
&+E_2(v)\nm_{S^c}(v)\nm(v)^{k_1+1}h^{(k_2+2)}(\rho,\nm(v)),
\end{align*}
with  $E_1(v)=(4+4k_1)w(v)+2(y_1\partial_{y_1}w(v)+z_1\partial_{z_1}w(v))$
and $E_2(v)=4w(v)$.  On dealing with the other terms in a similar fashion, this concludes the proof of the lemma when $\#S=1$.

Let $2\leq M\leq r_1+r_2$. 
For the inductive step let us assume the veracity of the lemma for all subsets of $\{1,\ldots,r_1+r_2\}$ of cardinality at most $M-1$.
Now let  $S$ be any subset of $\{1,\ldots,r_1+r_2\}$ with $M$ elements. We write
$S=S_1\cup \{j\}$ for
some $S_1$ of cardinality  $M-1$ elements. 
For simplicity we shall assume that $j\leq r_1$,  the complex case being handled similarly. 
The induction hypothesis implies that
there  
exist smooth weights
$w_{S_1,k_1,k_2}, 
w_{S_1,k_1,k_2}^{(m)}\in \cW_1(V)$, such that 
\begin{align*}
\partial_{S_1}K(v)
=~&
w_{S_1,k_1,k_2}(v)\nm(v)^{k_1}h^{(k_2)}(\rho,\nm(v))\\
& +
\sum_{m=0}^{d(S_1)}
w_{S_1,k_1,k_2}^{(m)}(v)
\nm_{S_1^c}(v)\nm(v)^{\max\{0,k_1+m-1\}}h^{(k_2+m)}
(\rho,\nm(v)).
\end{align*} 
Taking the derivative of this with respect to $v_j$, we obtain
\begin{align*}
\partial_{S} K(v)
=~&
\partial_j \left(w_{S_1,k_1,k_2}(v)\nm(v)^{k_1}h^{(k_2)}(\rho,\nm(v))\right)\\
& +
\sum_{m=0}^{d(S_1)}
\partial_j\left(w_{S_1,k_1,k_2}^{(m)}(v)
\nm_{S_1^c}(v)\nm(v)^{\max\{0,k_1+m-1\}}h^{(k_2+m)}
(\rho,\nm(v))\right).
\end{align*}
An application of the induction hypothesis shows that the first term here is satisfactory for the lemma. Turning to the $m$th summand, 
 we see that
\begin{align*}
\partial_j\Big(
&
w_{S_1,k_1,k_2}^{(m)}(v)
\nm_{S_1^c}(v)\nm(v)^{\max\{0,k_1+m-1\}}h^{(k_2+m)}
(\rho,\nm(v))\Big) \\
=~&
\partial_j(\nm_{S_1^c}(v))
\Big(
w_{S_1,k_1,k_2}^{(m)}(v)
\nm(v)^{\max\{0,k_1+m-1\}}h^{(k_2+m)}
(\rho,\nm(v))\Big) \\
&+
\nm_{S_1^c}(v)
\partial_j\Big(
w_{S_1,k_1,k_2}^{(m)}(v)
\nm(v)^{\max\{0,k_1+m-1\}}h^{(k_2+m)}
(\rho,\nm(v))\Big).
\end{align*}
Since  $j\leq r_1$, we have 
$\partial_j(\nm_{S_1^c}(v))=\nm_{S^c}(v)$ here.
Thus the first term is already in the desired form and we need to investigate the second term, $G(v)$, say.
Applying the 
induction hypothesis we obtain
weights 
$w_{S,k_1,k_2}^{(m)}, 
w_{S,k_1,k_2}^{(m,n)}\in \cW_1(V)$, such that 
\begin{align*}
G(v)=~&
\nm_{S_1^c}(v)
w_{S,k_1,k_2}^{(m)}(v) \nm(v)^{\max\{0,
k_1+m-1\}}h^{(k_2+m)}(\rho,\nm(v))\\
&+\nm_{S_1^c}(v)\sum_{n=0,1}\nm_{\{j\}^c}(v)
w_{S,k_1,k_2}^{(m,n)} (v)
\nm(v)^{\max\{0,k_1+m+n-2\}}h^{(k_2+m+n)}(\rho,\nm(v)).
\end{align*}
But $
\nm_{S_1^c}(v)
=
v_j\nm_{S^c}(v)$ and 
$
\nm_{S_1^c}(v)\nm_{\{j\}^c}(v)
=
\nm_{S^c}(v)\nm(v)$, whence
\begin{align*}
G(v)=~&
\nm_{S^c}(v)
v_jw_{S,k_1,k_2}^{(m)}(v) 
\nm(v)^{\max\{0,
k_1+m-1\}}h^{(k_2+m)}(\rho,\nm(v))\\
&+
\sum_{n=0,1}
\nm_{S^c}(v)
w_{S,k_1,k_2}^{(m,0)} (v)
\nm(v)^{\max\{1,
k_1+m+n-1\}}h^{(k_2+m+1)}(\rho,\nm(v)).
\end{align*}
This is satisfactory for the lemma and so concludes its proof.
\end{proof}

We shall ultimately be interested in a version of Lemma \ref{lemma1} in the special case $k_1=k_2=0$,  when $\partial_S$ is replaced by an arbitrary power of $\partial_S$.
This is an easy consequence of  our work so far, as the following result attests.

\begin{lemma}
\label{lemma1-cor}
Assume the same notation as in Lemma \ref{lemma1} and let $N\in \ZZ_{\geq 0}$.
Then there
exist weight functions 
$w_{S,N}, 
w_{S,N}^{(m,n)}\in \cW_1(V)$, such that 
\begin{align*}
\partial_{S}^N
\big\{w(v)h(\rho,\nm(v))\big\}
=~&
w_{S,N}(v)h(\rho,\nm(v))\\
& +
\sum_{n=1}^N
\sum_{m=0}^{Nd(S)}
\hspace{-0.2cm}
w_{S,N}^{(m,n)}(v)
\nm_{S^c}(v)^n \nm(v)^{\max\{0,m-n\}}h^{(m)}
(\rho,\nm(v)).
\end{align*} 
\end{lemma}

\begin{proof}
We argue by induction on $N\geq 1$, the case $N=1$ following from  Lemma \ref{lemma1}.
The induction hypothesis ensures that
\begin{align*}
\partial_{S}^{N+1}\big\{w(v)h(\rho,\nm(v))\big\}
=~&
\tilde{w}_{S,N}(v)
 +
\sum_{n=1}^N
\sum_{m=0}^{Nd(S)}
\nm_{S^c}(v)^n
\tilde{w}_{S,N}^{(m,n)}(v),
\end{align*} 
where
\begin{align*}
\tilde{w}_{S,N}(v)&=
\partial_S(w_{S,N}(v)h(\rho,\nm(v))),\\
\tilde{w}_{S,N}^{(m,n)}(v)&=
\partial_S\left(w_{S,N}^{(m,n)}(v)
 \nm(v)^{\max\{0,m-n\}}h^{(m)}
(\rho,\nm(v))\right),
\end{align*}
for 
$w_{S,N}, 
w_{S,N}^{(m,n)}\in \cW_1(V)$. 
Invoking Lemma \ref{lemma1} to 
evaluate $\tilde{w}_{S,N}$ and
$\tilde{w}_{S,N}^{(m,n)}$, 
it is a simple matter to check that one arrives at an expression suitable for the conclusion of the lemma.
\end{proof}

\subsection{Estimation of $p_\rho(v)$}

We now apply  our work in the previous section to the task of estimating $p_\rho(v)$, as 
given by \eqref{eq:p}. Let $\rho\ll 1$, with $\log \rho$ having order of magnitude $\log P$.
The following result shows that 
$p_\rho (v)$ is essentially supported on the set of $v\in V$ for which 
$\cH(v)\ll \rho^{-1}P^\ve$, in the notation of \eqref{eq:Tv'}.

\begin{lemma}
\label{pbsupport}
Let $\ve>0$ and $N\in \ZZ_{\geq 0}$.
Then we have 
$$
p_\rho(v)
\ll_{N}
\rho^{-1}\left(\rho^{-1}P^{\ve}|\cH(v)|^{-1}\right)^{N}.
$$
\end{lemma}
\begin{proof}
Recall that 
$$
p_\rho (v)=\int_{V} w_2(x)h(\rho,\nm(x))\e(-vx)\d x,
$$
where 
$w_2\in \cW_1(V)$. In particular  $\langle x\rangle\ll 1$ for all $x\in
\supp(w_2)$. This in turn implies that $\nm_S(x)\ll 1$ for any $x\in \supp(w_2)$ and any
subset $S$ of $\{1,\ldots,r_1+r_2\}.$

Let $S\subseteq \{1,\ldots,r_1+r_2\}$ and  let $N\in \ZZ_{\geq 0}$. Lemma \ref{lemma1-cor} implies that 
\begin{equation}\label{eq:linden-breakfast}
\partial_{S}^N \{w_2(x)h(\rho,\nm(x)\}
\ll_N
\sum_{n=1}^N 
\sum_{m=0}^{Nd(S)}
|\nm(x)|^{\max\{0,m-n\}}
|h^{(m)}(\rho,\nm(x))|.
\end{equation}
Integration by parts yields
\begin{equation}\label{eq:linden-lunch}
p_\rho (v)
\ll_N
|\nm_S(v)|^{-N}
\int_{\langle x \rangle \ll 1} \left|
\partial_S^N\left\{w_2(x)h(\rho,\nm(x)\right)\}\right|
\d x.
\end{equation}
Let $\ve>0$. We split the  integration here into two parts $J_1+J_2$, where 
$J_1$ arises from $|\nm(x)|\leq \rho P^{\ve/d}$ and 
$J_2$ is the contribution from  $|\nm(x)|>\rho P^{\ve/d}$.

Beginning with the latter, we deduce from  
taking $\nm(x)\ll 1$ in \eqref{eq:linden-breakfast} that 
$$
J_2\ll_N 
\sum_{m=0}^{Nd(S)}
\int_{\substack{
\langle x\rangle\ll
1\\
|\nm(x)|\geq
\rho P^{\ve/d}}} |h^{(m)}(\rho,\nm(x))| \d x.
$$
When $\rho\ll P^{\ve/d}$ the domain of integration is empty and so $J_2=0$.  Alternatively, 
if $\rho\gg P^{\ve/d}$ then 
Lemma \ref{lem:decay} implies that 
$h^{(m)}(\rho,\nm(x))\ll_{m,M} \rho^{-1-m}P^{-M{\ve/d}}, $ for any $M\in \ZZ_{\geq 0}$. 
Recall here that 
$\log \rho$ has order $\log P$. Thus, on taking $M$ sufficiently large, 
we have 
$
J_2\ll_{N}   P^{-N{\ve/d}}.
$
Turning to the estimation of $J_1$, 
Lemma \ref{lem:decay} yields $h^{(m)}(\rho, \nm(x))\ll \rho^{-m-1}$ for any $m\in \ZZ_{\geq 0}$. Hence \eqref{eq:linden-breakfast} gives
\begin{align*}
\partial_{S}^N \{w_2(x)h(\rho,\nm(x)\}
&\ll_N \sum_{n=0}^N \rho^{-1-n}\sum_{m=0}^{Nd(S)}(\rho^{-1}|\nm(x)|)^{\max\{0,m-n\}}\\
&\ll_N \rho^{-1-N}\sum_{m=0}^{Nd(S)}|\rho^{-1}\nm(x)|^{m}.
\end{align*}
We therefore obtain
\begin{align*}
J_1 
&\ll_{N}
\rho^{-1-N}
\sum_{m=0}^{Nd(S)}
\int_{\substack{
\langle x\rangle\ll 1\\
|\nm(x)|\leq \rho P^{\ve/d}}}
|\rho^{-1}\nm(x)|^{m}\d x\\
&\ll_{N}
\rho^{-1-N}
P^{N\ve},
\end{align*}
since $d(S)\leq r_1+r_2\leq d$.

Combining our estimates for $J_1,J_2$  in \eqref{eq:linden-lunch}, we therefore obtain
$$
p_\rho (v)
\ll_{N}
|\nm_S(v)|^{-N} \left( P^{-N\ve/d}+
\rho^{-1-N}
P^{N\ve}\right).
$$
Since $\rho\ll 1$ the second term here clearly dominates the first and we therefore conclude the proof of the lemma on taking $S=T(v)$, in the notation of \eqref{eq:Tv}  and 
\eqref{eq:Tv'}.
\end{proof}

Lemma \ref{pbsupport} will be effective when $\cH(v)$ is large, but we will need a companion ``trivial'' estimate to deal with the remaining cases. This is provided by the following result.

\begin{lemma}
We have $p_\rho(v)\ll |\log \rho|^{r_1+r_2-1}$. 
\label{pb bound}
\end{lemma}

\begin{proof}
We break the integral over $x=(x^{(1)}, \ldots, x^{(r_1+r_2)} )\in V$ 
in  \eqref{eq:p} 
into two parts $I_1+I_2$, say, where 
$I_1$ is the contribution from $x$ such that  $|x^{(l)}|< \rho$ for some $1\leq l \leq r_1+r_2$,  and 
$I_2$ has $|x^{(l)}|\geq \rho$ for every $1\leq l\leq r_1+r_2$.
Since $\supp(w_2)\ll 1$ we deduce from Lemma~\ref{lem:decay} with $i=j=N=0$ that
$I_1\ll 1.
$

Next, we take $i=j=0$ and   $N=2$ in Lemma \ref{lem:decay} to
 deduce that 
\begin{align*}
I_2
&\ll 1+ \rho^{-1} \int_{\substack{\langle x \rangle \ll 1\\
\min_l |x^{(l)}|\geq \rho}}  \min\left\{1, \frac{\rho}{|\nm(x)|}\right\}^2 \d x\\
&\ll 1+ \rho^{-1} \int_{\substack{\langle x \rangle \ll 1\\
\min_l |x^{(l)}|\geq \rho\\
\rho> |\nm(x)|
}}  \d x +
\rho \int_{\substack{\langle x \rangle \ll 1\\
\min_l |x^{(l)}|\geq \rho\\
\rho\leq |\nm(x)|
}} \frac{1}{|\nm(x)|^2} \d x\\
&\ll 1+ |\log  \rho|^{r_1+r_2-1}.
\end{align*}
This  completes the proof of the lemma since $|\log \rho|\gg 1$.
\end{proof}

\subsection{Application of Skinner's estimates}
In this section we show how Skinner's treatment of cubic exponential integrals 
in \cite[\S 6]{skinner} can be recycled here to help deal with the integral $K(v,\m)$ that appears in \eqref{ibc}.
Recall the definition \eqref{eq:omega} of the weight function $\omega$ and define
$$
I(v,\m)=\int_{V^n} \omega(\x/P)\e(vF(\x)-\m.\x)\d\x,
$$
for $v\in V$ and $\m \in V^n$.
Note that this is equal to the integral $I(v,-\m)$ introduced in \cite[Eq.\ (5.7)]{skinner}.
 Then it easily follows that 
$$ 
|I(v,\m)-K(v,\m)|\ll e^{-(\log P)^2/2}.
$$
Returning to \eqref{eq:no} and \eqref{ibc}, we shall in this section mainly be concerned with the 
contribution from non-zero phases
\begin{align*}
M(P) &=
\sum_{\substack{(0)\neq \mathfrak{b}\subseteq \fo\\
\n \mathfrak{b}  \ll Q^{d}
}} \sum_{
\substack{\0\neq \m\in \hat{\mathfrak{b}  }^n\\
\langle \m \rangle \leq P^A}}
(\n \mathfrak{b}  )^{-n}S_\mathfrak{b}  (\m)
I_\mathfrak{b}  (\m)\\
&=
\sum_{\substack{(0)\neq \mathfrak{b}\subseteq \fo\\
\n \mathfrak{b}  \ll Q^{d}
}} \sum_{
\substack{\0\neq \m\in \hat{\mathfrak{b}  }^n\\
\langle \m \rangle \leq P^A}}
(\n \mathfrak{b}  )^{-n}S_\mathfrak{b}  (\m)
\int_V p_\rho(v)
K(Q^{-2}v,\m) \d v.
\end{align*}
Now
$$
\int_V p_\rho(v)
K(Q^{-2}v,\m) \d v-
\int_V p_\rho(v)
I(Q^{-2}v,\m) \d v
\ll e^{-(\log P)^2/2} |\log \rho|^{r_1+r_2-1},
$$
by \eqref{pb bound}.
Since $\rho=Q^{-d}\n\fb \geq P^{-3d/2}$ we see that 
$|\log \rho|\ll \log P$. Observing that 
$ e^{-(\log P)^2/2} $ decays faster than any power of $P$, we may conclude that 
\begin{equation}\begin{split}
\label{i1zb}
M(P)&=
\sum_{\substack{(0)\neq \mathfrak{b}\subseteq \fo\\
\n \mathfrak{b}  \ll Q^{d}
}}
\int_V
 \sum_{
\substack{\0\neq \m\in \hat{\mathfrak{b}  }^n\\
\langle \m \rangle \leq P^A}}
(\n \mathfrak{b}  )^{-n}S_\mathfrak{b}  (\m)
 p_\rho(v)
I(Q^{-2}v,\m) \d v +O(1).
\end{split}
\end{equation}
Let $J^{(l)}(z,\m)$ denote the $l$th
component of $I(z,\m)$, so that 
$I(z,\m)=\prod_l J^{(l)}(z,\m)$.
We are now in a position to apply 
\cite[Eq.\ (6.32)]{skinner}, which gives the following result.

\begin{lemma}
Let $\ve>0$ and 
$B^{(l)}(z)=P^\ve(P^{-1}+|z^{(l)}|P^2 )$. Then we have 
$$
J^{(l)}(z,\m) \ll
\begin{cases}
		P^\ve\min\{P^{nc_l},|Pz^{(l)}|^{-nc_l/2} \},  & \mbox{if  $|\m^{(l)}|\ll
B^{(l)}(z)$},\\
		\exp(-\frac{1}{3}(\log P(2+|\m^{(l)}|))^2), & \mbox{otherwise}.
		\end{cases}
$$
In particular, when $|\m^{(l)}| \ll \Bl(z)$ for all $l$, we have
\begin{equation*}
I(z,\m)\ll P^{d\ve}\prod_{l=1}^{r_1+r_2}\min\{P^{nc_l},|P z^{(l)}|^{-nc_l/2}  \}.
\end{equation*}
\label{I bound}
\end{lemma}

Roughly speaking, the intuition behind the proof of this result is that one can use integration by parts when $|\m^{(l)}|$ dominates. Alternatively, when $|v^{(l)}\nabla F^{(l)}(\x^{(l)})|$ dominates methods from complex analysis are used to study the integral.
Lemma \ref{I bound} 
allows us to freely truncate the $\m$ summation in \eqref{i1zb} to $\m$ satisfying
\begin{equation}\label{eq:linden-dinner}
|\m^{(l)}|\ll B^{(l)}(Q^{-2}v) = P^{-1+\ve}(1+|v^{(l)}|),
\end{equation}
with acceptable error.  For such $\m$ we deduce that 
\begin{align*}
I(Q^{-2}v,\m)
&\ll P^{d\ve}\prod_l\min\{P^{nc_l},|PQ^{-2}v^{(l)}|^{-nc_l/2}  \}\\
&= P^{d\ve}\prod_l\min\{P^{nc_l},|P^{-2}v^{(l)}|^{-nc_l/2}  \}\\
&\ll P^{dn+d\ve}\cH(v)^{-n/2},
\end{align*} 
where  $\cH(v)$ is given by \eqref{eq:Tv'}.
Inserting this into \eqref{i1zb} yields
\begin{align*}
M(P) &\ll 1 +
P^{dn+d\ve}
\sum_{\substack{(0)\neq \mathfrak{b}\subseteq \fo\\
\n \mathfrak{b}  \ll Q^{d}
}} 
\int_V 
\cH(v)^{-n/2}
\sum_{
\substack{\0\neq \m\in \hat{\mathfrak{b}  }^n\\
\mbox{\scriptsize{\eqref{eq:linden-dinner} holds}}}}
(\n \mathfrak{b}  )^{-n}|S_\mathfrak{b}  (\m)|
|p_\rho(v)|
\d v.
\end{align*}
Let 
\[
R_\rho=
\{v\in V: \cH(v)\ll \rho^{-1}P^{2\ve}\}.
\]
Using 
Lemma \ref{pbsupport} to estimate $p_\rho(v)$,
we see that 
the overall contribution to the above estimate from 
$v\in V\setminus R_\rho$ is $O(1)$. 
When $v\in R_{\rho}$ we will simply invoke Lemma \ref{pb bound}, which gives
$p_\rho(v)\ll |\log \rho|^{r_1+r_2-1}\ll P^{\ve}$.
Hence we may write
\begin{align*}
M(P)&\ll 1+ 
P^{dn+(d+1)\ve}
\sum_{\substack{(0)\neq \mathfrak{b}\subseteq \fo\\
\n \mathfrak{b}  \ll Q^{d}
}} 
\int_{R_{\rho}} 
\cH(v)^{-n/2}
\sum_{
\substack{\0\neq \m\in \hat{\mathfrak{b}  }^n\\
\mbox{\scriptsize{\eqref{eq:linden-dinner} holds}}}}
(\n \mathfrak{b}  )^{-n}|S_\mathfrak{b}  (\m)|
\d v.
\end{align*}

Notice that $R_\rho\subseteq 
\{v\in V: \cH(v)\ll P^{3d/2+2\ve}\}=R$, say.
Let $\chi_\rho(v)$ be the characteristic function of $R_{\rho}$. 
Then we may write
\begin{equation}\label{eq:shot}
M(P)\ll 1+ 
P^{dn+(d+1)\ve}
\int_{R} E(v,P) \d v, 
\end{equation}
where
$$
E(v,P)=
\cH(v)^{-n/2}
\sum_{\substack{(0)\neq \mathfrak{b}\subseteq \fo\\
\n \mathfrak{b}  \ll Q^{d}
}} 
  \chi_\rho(v)
\sum_{
\substack{\0\neq \m\in \hat{\mathfrak{b}  }^n\\
\mbox{\scriptsize{\eqref{eq:linden-dinner} holds}}}}
(\n \mathfrak{b}  )^{-n}|S_\mathfrak{b}  (\m)|.
$$
We proceed to show that the outer sum  is actually restricted to
$Q_0\ll \n \fb \ll Q_1$, for suitable  $Q_0$ and $Q_1$.

Notice that  for $v\in R$, we have $v\in R_\rho$ if and only if 
$\cH(v)\ll \rho^{-1} P^{2\ve }$.
Recalling that  $\rho=Q^{-d}\n\bbb$,  we see that  $v\in R_\rho$ if and only if 
$$
\n\bbb\ll P^{2\ve} Q^d \cH(v)^{-1}=P^{3d/2+2\ve}\cH(v)^{-1}=Q_1,
$$ 
say.  
Recall 
from \eqref{eq:linden-dinner} that 
 $|\m^{(l)}|\ll P^{-1+\ve}(1+|v^{(l)}|)=B^{(l)}$, say,  for 
$\m=(m_1,\ldots,m_n)\in \hat{\fb}^n
$. 
Hence  
$
|\nm(m_i)|\ll 
|\nm(B)| \ll  P^{-d(1-\ve)}\cH(v),
$
for $1\leq i\leq n$.
Let 
$\m\in \hat{\fb}^n$ be  non- zero and suppose, without loss of generality, that  
$m_1\neq 0$. Then 
$m_1\in \hat\fb=\fb^{-1}\fd^{-1}$, and so 
$m_1 \fb\fd$ is integral. It follows that 
$\nm(m_1)\nm(\fb)\geq 1$ and so  $\nm(m_1)\geq \nm(\fb)^{-1}.$
This in turn  implies that
$$
\n\bbb\gg P^{d(1-\ve)}\cH(v)^{-1}=Q_0,
$$ 
say. Bringing this all together in \eqref{eq:no} and \eqref{eq:shot}, 
and replacing $2\ve$ by $\ve$, 
we may now record the following result.

\begin{lemma}
\label{lemma pb}
Let $\ve>0$ and let $B=\bigoplus_l B^{(l)}$, with  $B^{(l)}=P^{-1+\ve}(1+|v^{(l)}|) $. 
Let
$$
Q_0=P^{d(1-\ve)}\cH(v)^{-1}, \quad 
Q_1=
P^{3d/2+\ve}\cH(v)^{-1}.
$$
Then we have
$$
\No=
\frac{c_{Q}2^{r_2 n}}{D_{K}^{n/2}Q^{2d}}
\sum_{\substack{(0)\neq \mathfrak{b}\subseteq \fo\\
\n \mathfrak{b}  \ll Q^{d}
}} (\n \mathfrak{b}  )^{-n}S_\mathfrak{b}  (\0)
I_\mathfrak{b}  (\0) 
+
O\left(1+ P^{d(n-3)+(d+1)\ve}
\int_{R} E(v,P) \d v\right),
$$
with $R=\{v\in V: \cH(v)\ll P^{3d/2+\ve}\}$ and 
$$
E(v,P)=
\cH(v)^{-n/2}
\sum_{\substack{(0)\neq \mathfrak{b}\subseteq \fo\\
Q_0\ll \n \mathfrak{b}  \ll Q_1
}} 
\sum_{
\substack{\0\neq \m\in \hat{\mathfrak{b}  }^n\\
|\m^{(l)}|\ll B^{(l)}}}
(\n \mathfrak{b}  )^{-n}|S_\mathfrak{b}  (\m)|.
$$
\end{lemma}

In applying this result the goal is to show that 
$\int_R E(v,P)\d v\ll  P^{-\Delta}$,  for a suitable constant $\Delta>0$.
When $\cH(v)$ is small, one sees that $Q_0$ is large and one can hope to gain sufficient cancellation in the  exponential sums $S_\fb(\m)$.
On the other hand, when $\cH(v)$ is large, then $Q_1$ is small and the 
exponential sums involved are small. However,  in this case,  the factor 
$\cH(v)^{-n/2}$ will produce the necessary saving.  This is the main idea behind our application of 
Lemma~\ref{lemma pb}.

\section{Cubic exponential sums}\label{s:sums}

The purpose of this section is to make a careful analysis of the exponential sums 
$S_\fb(\m)$ 
occurring in Theorem  \ref{prop:2} and Lemma \ref{lemma pb}
 when $F\in \fo[X_1,\ldots,X_n]$ is a
non-singular cubic form, with $n\geq 3$. Here $\fb$ is an arbitrary integral ideal and $\m\in {\hat
\fb}^n$.
Applying 
Lemma \ref{lem:orthogonal}, we see that there 
exists $\gamma=\nu/\alpha$, 
 for $\nu,\alpha\in \fo$ such that $(\nu)$ is coprime to $\fb$, which allows us to write
\begin{align*}
S_\mathfrak{b}  (\m)
&=
\sum_{a\in (\fo/\fb)^*} 
\sum_{\ma{a}\bmod{\mathfrak{b} }}\e\left( a\gamma F(\ma{a})+ \m. \ma{a}\right).
\end{align*}
It follows from our work in \S \ref{s:characters} that  
\begin{equation}\label{eq:construct}
\fa_\gamma=\fb\fd.
\end{equation}
Moreover, there exists a prime ideal $\fp_1$ coprime to $\fb\fd$, such that $(\alpha)=\fb\fd\fp_1$.
The following standard result, established by Skinner \cite[Lemma 3]{skinner}, will prove useful in
our analysis. 

\begin{lemma}\label{lem:skinner-3}
Suppose that $\fa,\fb,\fc$ are integral ideals such that $\fa=\fb\fc$. Then we have the following:
\begin{enumerate}
\item[(i)]
if $\alpha\in \fo$ satisfies $\ord_\fp(\alpha)=\ord_\fp(\fb)$ for all $\fp\mid \fa$, then 
$$
\fo/\fa=\{\beta+\mu \alpha:  \beta\in \fo/\fb, ~\mu\in \fo/\fc\}\mbox{;}
$$
\item[(ii)]
if, furthermore,  $\fb$ and $\fc$ are coprime and 
if $\alpha, \lambda\in \fo$ satisfy $\ord_\fp(\alpha)=\ord_\fp(\fb)$ 
and 
$\ord_\fp(\lambda)=\ord_\fp(\fc)$ 
for all $\fp\mid \fa$, then 
$$
\fo/\fa=\{\alpha \mu +\lambda \beta:  \beta\in \fo/\fb, ~\mu\in \fo/\fc\}.
$$
\end{enumerate}
\end{lemma}

There is an abuse of notation at play in this lemma, in that the sets involved are actually coset
representatives for $\fo/\fa, \fo/\fb$ and $\fo/\fc$.

Since we aim to monopolise upon the existing work of Skinner,  
we will begin by  indicating how our expression for $S_{\fb}(\m)$ is related to the  exponential
sums 
$$
S(\gamma,\b) = 
\sum_{\ma{a}\bmod{\mathfrak{a}_\gamma }}\e\left( \gamma F(\ma{a})- \b. \ma{a}\right),
$$
that emerge in  \cite[Eq.\ (5.6)]{skinner}, 
for $\b\in \hat \fa_\gamma$.
This is the object of  the following result.

\begin{lemma}\label{lem:ram}
For any integral ideal $\fb$ and $\m\in {\hat \fb}^n$
we have 
$$
S_{\fb}(\m)
=
D_K^{-n}
\sum_{a\in (\fo/\fb)^*}  S(a\gamma,-\m).
$$
\end{lemma}

\begin{proof}
Recall from \eqref{eq:construct} that
$\fa_\gamma=\fb\fd$.
Hence for any $\m\in \hat{\bbb}^n$ we have 
\begin{align*}
\sum_{a\in (\fo/\fb)^*}  S(a\gamma,-\m)
&=
\sum_{a\in (\fo/\fb)^*} 
\sum_{\a\bmod{\fb\fd}}\e(a\gamma F(\a)+\m.\a).
\end{align*}
By Lemma \ref{lem:alg1} (i) we may find 
$\beta\in \fb$ such that $\ord_\fp(\beta)=\ord_\fp(\fb)$, for each prime ideal
$\fp\mid \fb\fd$. Then Lemma \ref{lem:skinner-3}
allows us to write 
$\a=\b+\beta\c$, with $\b$ running modulo $\fb$ and $\c$ running modulo $\fd$,
 giving 
\begin{align*}
\sum_{a\in (\fo/\fb)^*}  S(a\gamma,-\m)
&=
\sum_{a\in (\fo/\fb)^*}  
\sum_{\c\bmod{\fd}}\sum_{\b\bmod{\fb}}\e(a\gamma F(\b+\beta\c)+\m.(\b+\beta\c)).
\end{align*}
Since $\fa_\gamma=\fb\fd$ and $\fb\mid (\beta)$, it follows that 
$\gamma \beta\in \hat\fo$,
whence
$\e(a\gamma F(\b+\beta\c))=\e(a\gamma F(\b)).$ 
Moreover, if 
 $\m\in \hat{\bbb}^n$ then one  has
$\e(\beta\m.\c)=1$,  since
$ \beta \hat{\fb}=\beta\fb^{-1}\hat\fo\subseteq \hat\fo$. 
Thus
\begin{align*}
\sum_{a\in (\fo/\fb)^*}  S(a\gamma,-\m)
&= (\n \fd)^{n}
\sum_{a\in (\fo/\fb)^*}  
\sum_{\b\bmod{\fb}}\e(a\gamma F(\b)+\m.\b)\\
&= D_K^{n}
S_{\fb}(\m),
\end{align*}
as required.
\end{proof}

It will be convenient to pass from exponential sums modulo $\fb$  indexed by $\hat\fb^n$, to exponential sums modulo $\fb$ indexed by $\fo^n$.
Define
\begin{equation}\label{eq:tilde-S}
\tilde S_\mathfrak{b}  (\v)=
\sum_{a\in (\fo/\fb)^*} 
\sum_{\ma{a}\bmod{\mathfrak{b}}}\e\left( \gamma\{ a F(\ma{a})+ 
 \v. \ma{a}\}\right),
\end{equation}
for any $\v\in \fo^n$.
It follows from our work in \S \ref{s:characters} that this expression is independent of the precise
choice of $\gamma$. 
Since $(\nu)$ is coprime to $\fb$, we may write
\begin{equation}\label{eq:call}
\begin{split}
S_\mathfrak{b}  (\m)
&=
\sum_{a\in (\fo/\fb)^*} 
\sum_{\ma{a}\bmod{\mathfrak{b}}}
\e\left( a\gamma F(\nu \ma{a})+ \nu \m. \ma{a}\right)
= \tilde S_\mathfrak{b}  (\alpha \m),
\end{split}
\end{equation}
on making a suitable change of variables in the $a$ and $\a$ summations.
Note here that $\alpha \m$ belongs to $\fo^n$. In fact, following the notation introduced at the start of this section, one easily sees that $\alpha\m$ belongs to $\fp_1^n$.

The exponential sums \eqref{eq:tilde-S} satisfy the basic multiplicativity property
\begin{equation}\label{eq:Lucy}
\tilde S_{\fb_1\fb_2}  (\v)=\tilde S_{\mathfrak{b}_1}  (\v)\tilde S_{\mathfrak{b}_2}  (\v),
\end{equation}
if $\fb_1,\fb_2$ are coprime integral ideals. We will use this decomposition to analyse $\tilde S_\fb(\v)$
at the square-free and square-full parts of $\fb$ separately.
To check this we write $\fb=\fb_1\fb_2$ and recall that $\fa_\gamma=\fb \fd$.
Applying Lemma \ref{lem:alg1}(i)
we find   $\alpha $ and $\lambda$ such that $\ord_{\fp}(\alpha)=\ord_\fp(\fb_1)$ and $\ord_{\fp}(\lambda)=\ord_\fp(\fb_2)$ for all $\fp \mid \fb\fd $.  Lemma \ref{lem:skinner-3}(ii) then allows us to write $a=\lambda b+\alpha c$ and $\a=\lambda\b+\alpha\c$, for $b,\b \bmod{\fb_1}$
and $c,\c \bmod{\fb_2}$,  to get
\begin{align*}
\tilde S_\fb(\v)
&=\sum_{b\in (\fo/\fb_1)^*} \sum_{c\in (\fo/\fb_2)^*}
\sum_{\b\bmod{\mathfrak{b}_1}}\sum_{\c\bmod{\mathfrak{b}_2}}
\e\left( \gamma \{(\lambda b+\alpha c)F(\lambda\b+\alpha\c)+ \v.(\lambda\b+\alpha\c)\}\right)\\
&=\sum_{b\in (\fo/\fb_1)^*}\sum_{\b\bmod{\mathfrak{b}_1}}\e\left(\lambda\gamma
\{b\lambda^3F(\b)+\v.\b\}\right) \sum_{c\in (\fo/\fb_2)^*}
\sum_{\c\bmod{\mathfrak{b}_2}}
\e\left(\alpha \gamma \{c\alpha^3F(\c)+ \v.\c\}\right)\\
&=\sum_{b\in (\fo/\fb_1)^*}\sum_{\b\bmod{\mathfrak{b}_1}}\e\left(\lambda\gamma
\{bF(\b)+\v.\b\}\right) \sum_{c\in (\fo/\fb_2)^*}
\sum_{\c\bmod{\mathfrak{b}_2}}
\e\left(\alpha \gamma \{cF(\c)+ \v.\c\}\right).
\end{align*}
The multiplicativity property therefore follows on noting that 
$\fa_{\lambda\gamma}=\fb_1\fd $ and $\fa_{\alpha\gamma}=\fb_2\fd$.

\subsection{Square-free $\fb$}

Let $G\in \fo[X_1,\dots,X_n]$ be the dual form of $F$, whose zero locus parameterises the set of
hyperplanes whose intersection with the cubic hypersurface $F=0$ produce a singular variety. It is
well-known that $G$ is absolutely irreducible and has degree $3\cdot 2^{n-2}$. 
The primary aim of this section is to establish the following result, which is an exact 
analogue of \cite[Lemma 13]{hb-10}.

\begin{lemma}\label{lem:hb-13}
Let $\fb$ be a square-free integral ideal
and let $\v\in \fo^n$.
Then there exists an absolute constant $A>0$ such that 
$$
 |\tilde S_\bbb(\v)|\leq A^{\omega(\bbb)} (\n\bbb)^{(n+1)/2} (\n \fh)^{1/2},
$$
where
$\fh$ is the greatest common divisor of $\fb$ and $(G(\v))$.
\end{lemma}

Now for any square-free ideal $\fb$ it follows from \cite[Lemma 23]{skinner}, Lemma \ref{lem:ram} and \eqref{eq:call}
that
$$
 |\tilde S_\bbb(\v)|\leq A^{\omega(\bbb)} (\n\bbb)^{n/2+1}.
 $$
Hence, 
by multiplicativity, in order to complete the proof of Lemma \ref{lem:hb-13}
it will suffice to show that
\begin{equation}\label{eq:suffice}
\tilde S_\fp(\v)\ll (\n\fp)^{(n+1)/2},
\end{equation}
when $\fp$ is a prime ideal such that $\fp\nmid G(\v)$. 
We may further assume that $\fp$ is unramified since otherwise one has the trivial bound $\tilde S_{\fp}(\v)=O(1)$.  Our hypotheses on $\fp$ imply 
that $F$  and $F_{\v}$ are both non-singular modulo $\fp$, where 
$F_{\v}$ is the cubic form in $n-1$ variables obtained by eliminating a variable from 
the pair of equations $F(\mathbf{X})=0$ and $\v . \mathbf{X}=0$.
In order to bound $\tilde S_\fp(\v)$ we  introduce a dummy sum over an extra variable $t$ to get
\begin{align*}
\phi(\fp)\tilde S_{\fp}(\v)
&=
\sum_{a\in (\fo/\fp)^*} 
\sum_{t\in (\fo/\fp)^*} 
\sum_{\a\bmod{\mathfrak{p}}} \e(\gamma\{ a\overline{t}^3 F(\a)+ \v. \a\})\\
&=
\sum_{a\in (\fo/\fp)^*} 
\sum_{t\in (\fo/\fp)^*} 
\sum_{\a\bmod{\mathfrak{p}}}
\e(\gamma\{ aF(\a)+t\v. \a\})\\
&=
\sum_{a\in (\fo/\fp)^*}
\left( 
\sum_{t\bmod{\fp}} 
\sum_{\a\bmod{\mathfrak{p}}}
\e(\gamma\{a F(\a)+t\v. \a\})-
\sum_{\a\bmod{\mathfrak{p}}}\e(a\gamma F(\a))\right),
\end{align*}
where $\overline t$ is the multiplicative inverse of $t$ modulo $\fp$.
But clearly 
\begin{align*}
\sum_{a\in (\fo/\fp)^*}
\sum_{t\bmod{\fp}} 
\sum_{\a\bmod{\mathfrak{p}}}
\e\left(\gamma \{ aF(\a)+t\v. \a\}\right)
&=
\n \fp \sum_{a\in (\fo/\fp)^*}
\sum_{\substack{\a\bmod{\mathfrak{p}}\\
\fp \mid \v.\a}}
\e(a \gamma F(\a)),
\end{align*}
since  $\v\in \fo^n$.
Assuming without loss of generality that $v_n\neq 0$, we may 
eliminate $a_n$ from this exponential sum, leading to  the expression
\begin{align*}
\phi(\fp)\tilde S_{\fp}(\v)=
\sum_{a\in (\fo/\fp)^*}
\left( 
\n \fp 
\sum_{\substack{\a' \bmod{\mathfrak{p}}}}
\e(a\gamma F_{\v}(\a'))
-\sum_{\a\bmod{\mathfrak{p}}}\e(a\gamma F(\a))\right),
\end{align*}
where $\a'=(a_1,\ldots,a_{n-1})$.  Since $F$ and $F_{\v}$ are non-singular modulo $\fp$, it
follows from 
Deligne's estimate \cite[Thm.~8.4]{deligne} that the 
two terms in the brackets are $O( (\n\fp)^{(n+1)/2})$. This therefore establishes
\eqref{eq:suffice}, which  concludes the proof of Lemma \ref{lem:hb-13}.

\subsection{Square-full $\fb$}

In this section we examine the exponential sum $\tilde S_{\bbb}(\v)$ in \eqref{eq:tilde-S}
for square-full integral ideals $\fb$ and suitable $\v\in \fo^n$,
the main idea being to average over the $\v$.
We begin with the following result.

\begin{lemma}\label{lem:skinner-21}
Let $B=\bigoplus_l B^{(l)} \in V$,  with $B^{(l)}> 0$. Let $\ve>0$ and let $\bbb$ be any square-full 
integral ideal. Let $\fc$ be an integral ideal which is coprime to $\fb\fd$, with $\fb\fd\fc$ principal.
Then we have 
$$
\sum_{\substack{\v\in \fc^n \\ 
|\v^{(l)}| \leq B^{(l)}}} |\tilde S_\fb(\v)| \ll
(\n\bbb)^{n/2+1+\ve}
\left(\nm(B)^n+(\n\fb)^{n/3}\right). 
$$ 
The implied constant in this estimate does not depend on $\fc$.
\end{lemma}

\begin{proof}
By hypothesis there exists $\alpha\in \fo$ such that 
$\fb\fd\fc=(\alpha)$. Therefore
$\alpha\hat{\fb}=\fb\fd\fc(\fb\fd)^{-1}=\fc$ and it follows that any 
$\v\in \fc^n$ can be written  $\v=\alpha \m$ for $\m\in \hat \fb^n$.
Following our conventions, we may view $\alpha=(\alpha^{(1)},\ldots,\alpha^{(r_1+r_2)})$ as an 
element of $V$ by setting  $\alpha^{(j)}=\rho_j(\alpha)$. 
Hence \eqref{eq:call} and 
Lemma \ref{lem:ram} yield
\begin{align*}
\sum_{\substack{\v\in \fc^n \\ 
| \v^{(l)}|\leq B^{(l)}}} |\tilde S_\fb(\v)| 
&=\sum_{\substack{\v\in \fc^n \\ 
| \v^{(l)}|\leq B^{(l)}}} | S_\fb(\alpha^{-1}\v)| \\
&\leq 
\sum_{\substack{\m\in \hat\fb^n \\ 
| \m^{(l)} |\leq |\alpha^{(l)}|^{-1}B^{(l)}}} |S_\fb(\m)|\\
&\leq
\sum_{a\in (\fo/\fb)^*}
\sum_{\substack{\m\in \hat\fb^n \\ 
| \m^{(l)} |\leq |\alpha^{(l)}|^{-1}B^{(l)}}} 
  |S(a\gamma,-\m)|.
\end{align*}
Now it is clear that $\ord_{\fp}(\gamma)=\ord_{\fp}(a\gamma)$ for all $\fp\mid \fb$, since $a$ is coprime to $\fb$.
Thus  $\fb\mid \fa_{a\gamma}$  and it follows that $\hat\fb\subseteq \hat\fa_{a\gamma}$.
Enlarging the sum over $\m$, 
\cite[Lemma 21]{skinner} yields
\begin{align*}
\sum_{\substack{\v\in \fc^n \\ 
| \v^{(l)}|\leq B^{(l)}}} |\tilde S_\fb(\v)| 
&\ll
\sum_{a\in (\fo/\fb)^*}
(\n\fb)^{n/2+\ve}\left(\nm(\alpha)^{-n}\nm(B)^n (\n\fb)^n+(\n\fb)^{n/3}\right),
\end{align*}
for any $\ve>0$.  The lemma follows on noting that  
$\nm(\alpha)^{-1}\n\fb=(\n\fd\n\fc)^{-1}\leq 1$.
\end{proof}

We will need a companion estimate which deals with the problem of averaging over those $\v$ at which 
the dual form $G$ vanishes. The rest of this section will be devoted to a proof of the following
result. 

\begin{lemma}
\label{lem:hb-16}
Let $B=\bigoplus_l B^{(l)} \in V$,  with $B^{(l)}>0$. Let $\ve>0$ and let $\bbb$ be any square-full integral
ideal. Then we have 
\[
\sum_{\substack{\v\in \oh^n \\ 
| \v^{(l)}| \leq B^{(l)}\\
G(\v)=0
}} |\tilde S_\fb(\v)| \ll
(\n\fb \nm(B))^\ve
\left((\n\bbb)^{n+1/2}+
(\nm(B))^{n-3/2}(\n\bbb)^{n/2+4/3}\right).
\]
\end{lemma}

Lemmas \ref{lem:skinner-21} and \ref{lem:hb-16} are precise analogues of Lemmas 14 and 16,
respectively,  
of Heath-Brown \cite{hb-10}. The proofs for general number fields are very similar to the case
$K=\QQ$
and we shall attempt to be brief in our demonstration of Lemma \ref{lem:hb-16}.
The rationale behind this result is  the need to make up for the loss in Lemma
\ref{lem:hb-13} when $G( \v)$ vanishes.  Following Heath-Brown \cite[\S 7]{hb-10} we will get
extra savings from two sources:  firstly by summing non-trivially over $a$ and secondly by using the
relative sparsity of vectors $\v$ such that $G(\v)=0$.

Since $\fb$ is square-full 
we may  write  
$\fb=\fb_1^2\fb_2$, where $\fb_2$ is square-free and 
$\fb_2\mid \fb_1$.  Let us set $q=\n\fb$ and  $q_i=\n\fb_i$ for $i=1,2$.
Let $\gamma$ be as before and recall \eqref{eq:construct}.
Applying Lemma~\ref{lem:alg1}(i), 
let $\beta\in \fb_1$
and $\mu\in \fb_2$  such that $\ord_\fp(\beta)=\ord_\fp(\bbb_1)$ and 
$\ord_\fp(\mu)=\ord_\fp(\bbb_2)$, for all $\fp\mid \bbb\fd$. 
Notice that $\fa_{\gamma\beta^2\mu}=\fd$.
  Adapting the argument leading to   \cite[Eq.\ (8.5)]{skinner}, we obtain
\begin{align*}
 \tilde S_\bbb(\v)&=
 \sum_{a\in (\fo/\fb)^*}
 \sum_{\a\bmod \fb}\e(\gamma\{aF(\a)+\v.\a\} )\\
&=
\sum_{a\in (\fo/\fb)^*}
\sum_{\f\bmod{\bbb_1\bbb_2}}\sum_{\g \bmod{\bbb_1}}
\e(\gamma\{aF(\f+\beta\mu
\g)+\v.(\f+\beta\mu \g)\} )\\
&=
\sum_{a\in (\fo/\fb)^*}
\sum_{\f\bmod{\bbb_1\bbb_2}}
\e(\gamma\{aF(\f)+\v.\f)\} )
\sum_{\g \bmod{\bbb_1}}
\e(\gamma \beta \mu\{a
\g.\nabla F(\f)+\v.\g\}).
\end{align*}
Here $\e(\gamma \cdot )$
is  a primitive character modulo $\bbb$ and 
 $\e(\gamma\beta\mu \cdot )$
is  a primitive character modulo $\bbb_1$, since
 $\fa_{\gamma\beta\mu}=\fb_1\fd$.
 This implies that 
\begin{align*}
 \tilde S_\bbb(\v)&=
q_1^n
\sum_{a\in (\fo/\fb)^*}
\sum_{\substack{\f\bmod{\bbb_1\bbb_2}\\
a\nabla F(\f)+\v\in
\fb_1^n}}
\e(\gamma\{aF(\f)+\v.\f)\} ).
\end{align*}
Writing  $a=t+u\beta$, we obtain
\begin{align*}
 \tilde S_\bbb(\v)&= q_1^n
\sum_{t\in (\fo/\fb_1)^*} 
\sum_{\substack{\f\bmod{\fb_1\fb_2}\\ t\nabla F(\f)+\v\in \fb_1^n}}
e(\gamma\{tF(\f)+\v.\f\})\sum_{u\bmod{\fb_1\fb_2}}
\e(\gamma\beta uF(\f)) \\
&=q_1^{n+1}q_2
\sum_{t\in (\fo/\fb_1)^*} 
\sideset{}{^{(1)}}\sum_{\f}\e(\gamma\{ tF(\f)+\v.\f\}),
\end{align*}
where $\sum^{(1)}$ indicates that $\f$ runs modulo $\fb_1\fb_2$ 
subject to the constraints $t\nabla F(\f)+\v\in \fb_1^n $ and 
$F(\f)\in \fb_1\fb_2$. 
Here we have used the fact that $\e(\gamma \beta\cdot)$ is a primitive character modulo
$\fb_1\fb_2$.
Since $\beta\in \fb_1$ and $\fb_2\mid \fb_1$, it follows that $\beta^2\in \fb_1\fb_2$.
We next substitute $\f=\h+\beta \j$. Then, if $F(\h)=\beta m $, we will have $F(\f)\in \fb_1\fb_2$
precisely when $m+\j.\nabla F(\h)\in \fb_2 $. 
Moreover, if 
$t\nabla F(\h)=-\v+\beta \k $, then 
\[
tF(\f)+\v.\f\equiv tF(\h)+\v.\h+\beta^2(\k.\j+t\h.\nabla F(\j))\bmod{\fb}.
\]
It follows that 
\begin{align*}
\left|
\sideset{}{^{(1)}}\sum_{\f}
\e(\gamma\{ tF(\f)+\v.\f\})\right|
\leq
\sideset{}{^{(2)}}\sum_{\h}\max_{m,\k\bmod{\fb_2}}\left|\sideset{}{^{(3)}}\sum_{\j}\e(\gamma\beta^2
\{\k.\j+t\h.\nabla F(\j)\}) \right|,
\end{align*}
where $\sum^{(2)}$ is for $\h$ modulo $\fb_1$ such that 
$t\nabla F(\h)+\v\in \fb_1^n$ and $F(\h)\in \fb_1$,  and $\sum^{(3)}$ is over $\j$ modulo $\fb_2$ for
which $m+\j.\nabla F(\h)\in \fb_2$.

Let $\chi=\gamma\beta^2$
and note that $\fb_2\fd=\fa_{\chi}$.
Therefore $\e(\chi \cdot)$ denotes a primitive character modulo ${\fb_2}$.
We proceed by  bounding the inner sum over $\j$, which we write as $S^{(3)}$. 
By orthogonality we have
\begin{align*}
S^{(3)}
 &= q_2^{-1}\sum_{l\bmod{\fb_2}}\sum_{\j \bmod{\fb_2}}
 \e\left(\chi \{lm+(\k+l\nabla F(\h)).\j)+t\h.\nabla F(\j)\}\right)\\
 &\ll \max_{\mathbf{l} \bmod{\fb_2}}\left|
 \sum_{\j\bmod{\fb_2}}
 \e\left(\chi\{\mathbf{l}.\j+t\h.\nabla F(\j)\}\right)\right|.
\end{align*}
We proceed by adapting the  argument leading to \cite[Eq.\ (8.13)]{skinner}. 
Let $S_{\mathbf{l},\h}$ denote the sum over $\j\bmod{\fb_2}$.
Then
$$
|S_{\mathbf{l},\h}|^2=\sum_{\j_1\bmod{\fb_2}}\sum_{\j_2\bmod{\fb_2}}\e\left(
\chi\{\mathbf{l}.(\j_1-\j_2)+t\h.(\nabla
F(\j_1)-\nabla F(\j_2))\}\right).
$$
Write $\j_1=\j_2+\j_3$ to get
$$
\h.(\nabla F(\j_1)-\nabla F(\j_2))\equiv \h.\nabla F(\j_3)+6\j_2.\B(\h,\j_3)\bmod{\bbb_2},
$$
where $\mathbf{B}(\x,\y)$ is the system of $n$ bilinear forms defined in
\cite[\S 2]{skinner}.
It follows that
\begin{align*}
|S_{\mathbf{l},\h}|^2
&= 
\sum_{\j_3\bmod{\fb_2}}\e\left(\chi\{\mathbf{l}.\j_3+t\h.\nabla F(\j_3)\}\right)
\sum_{\j_2\bmod{\fb_2}}\e\left(6\chi t
\j_2.\B(\h,\j_3)\right) \\
&\ll q_2^n\# \{ \j_3\bmod{\bbb_2}:\B(\h,\j_3)\in (6^{-1}\bbb_2\cap\oh)^n\}.
\end{align*}
Thus we have shown that 
$S^{(3)}\ll q_2^{n/2}K(\fbth;\h)^{1/2}$,
where
$\fb_3=6^{-1}\fb_2\cap \fo$ and 
$$
K(\fb_3;\x)=\#\left\{\y\bmod{\fb_3}:\mathbf{B}(\x,\y)\in \fb_3^n\right\}.
$$ 
Once inserted into our work so far we therefore obtain
\begin{equation*}
 \tilde S_{\fb}(\v)\ll q_1^{n+1}q_2^{n/2+1}
 \sum_{t\in (\fo/\fb_1)^*}
 \sumtwo_{\h}K(\fb_3;\h)^{1/2}.
\end{equation*}

It is now time to introduce the summation over $\v$, defining
\begin{equation}\label{eq:def-N}
\mathcal{N}=\max_{\r\bmod{\fb_1}}\#\left\{\v\in \oh^{n}:|\v^{(l)}|\leq B^{(l)}, 
~G(\v)=0, ~\v\equiv \r\bmod{\fb_1}\right\}.
\end{equation}
Then we have 
\begin{align*}
\sum_{\substack{\v\in \oh^n \\ 
| \v^{(l)}|\leq B^{(l)}\\
G(\v)=0
}} |\tilde S_\fb(\v)| \ll
q_1^{n+2}q_2^{n/2+1}\mathcal{N}\sum\limits_{\substack{\h\bmod{\fb_1} \\ F(\h)\in
\fb_1}}K(\fb_3;\h)^{1/2}
=
q_1^{n+2}q_2^{n/2+1}\mathcal{N} S(\fb),
\end{align*}
say. The following result is concerned with an upper bound for $S(\fb)$.

\begin{lemma}\label{lem:S-c}
There exists an absolute constant $A>0$ such that 
$S(\fb)\leq   A^{\omega(\fb)} q_1^{n-1}q_2^{1/2}$.
\end{lemma}

The proof of Lemma \ref{lem:S-c} is the exact analogue of the treatment of $S(q)$ in \cite[\S
7]{hb-10}
and we have decided to omit the proof. Note that the analogue of \cite[Lemma 4]{hb-10}
(resp.\ the estimate $S_0(p^g;\0)\ll p^{(5n/6+1+\ve)g}$)
is provided by \cite[Lemma 10]{skinner} (resp.\ by Lemmas \ref{lem:hb-13} and~\ref{lem:skinner-21}).
Applying this result in our work above now yields
 $$
\sum_{\substack{\v\in \oh^n \\ 
| \v^{(l)}|\leq B^{(l)}\\
G(\v)=0
}} |\tilde S_\fb(\v)| \ll
A^{\omega(\fb)}q_1^{2n+1}q_2^{(n+3)/2} \mathcal{N},
$$
where $\mathcal{N}$ is given by \eqref{eq:def-N}. 
Finally, it follows from  Lemma \ref{lem:cohen} that
\begin{align*}
\sum_{\substack{\v\in \oh^n \\ 
| \v^{(l)}|\leq B^{(l)}\\
G(\v)=0
}} |\tilde S_\fb(\v)| \ll
(q \nm(B))^\ve q_1^{2n+1}q_2^{(n+3)/2}\left(1+\frac{\nm(B)}{q_1}\right)^{n-3/2},
\end{align*}
for any $\ve>0$.
Since $q=q_1^2q_2$ and $q_2\leq q^{1/3}$, the above bound is
 \begin{align*}
 &\ll(q \nm(B))^\ve\left\{
q^{n+1/2}+(\nm(B))^{n-3/2}(q_1^2q_2)^{n/2+5/4}q_2^{1/4}\right\}\\
&\ll (q \nm(B))^\ve\left\{ q^{n+1/2}+(\nm(B))^{n-3/2}q^{n/2+4/3}\right\}.
\end{align*}
This therefore concludes the proof of Lemma \ref{lem:hb-16}.

\section{Final deduction  of Theorem \ref{t:10}}\label{s:conclusion}

In this section we complete the proofs of Theorems \ref{t:10} and \ref{t:10'}, by combining our analysis of the exponential sums in \S \ref{s:sums} with Lemma \ref{lemma pb}.
In what follows it will be notationally convenient to follow the convention that the small positive  constant $\ve$ takes different values at different parts of the argument.

Recall from Lemma \ref{lemma pb} that 
$$
E(v,P)=
\cH(v)^{-n/2}
\sum_{\substack{(0)\neq \mathfrak{b}\subseteq \fo\\
Q_0\ll \n \mathfrak{b}  \ll Q_1
}} 
\sum_{
\substack{\0\neq \m\in \hat{\mathfrak{b}  }^n\\
|\m^{(l)}|\ll B^{(l)}}}
(\n \mathfrak{b}  )^{-n}|S_\mathfrak{b}  (\m)|.
$$
where $\cH(v)$ is given by \eqref{eq:Tv'}
and $B^{(l)}=P^{-1+\ve}(1+|v^{(l)}|)$, with  
$$
Q_0=P^{d(1-\ve)}\cH(v)^{-1}, \quad 
Q_1=
P^{3d/2+\ve}\cH(v)^{-1}.
$$
We will need to show that there exists an absolute constant $\Delta>0$, which is independent of $\ve$,  such that 
\begin{equation}\label{eq:white}
\int_R E(v,P) \d v \ll P^{-\Delta},
\end{equation}
where  $R=\{v\in V: \cH(v)\ll P^{3d/2+\ve}\}$.

We now write $\bbb=\bbb_1\bbb_2$,  where $\bbb_1$ is square-free and $\bbb_2$ is square-full. 
By Lemma \ref{lem:alg1}(ii) there exists $\alpha'\in \fb_2\fd$ and an unramified prime ideal $\fp_2$ 
coprime to $\fb_2$, such that $(\alpha')=\fb_2\fd\fp_2$ and $\n\fp_2\ll (\n \fb)^\ve$. 
Likewise, a second application of this result shows that there exists $\alpha\in \fb\fd$ and an unramified prime ideal $\fp_1$ 
coprime to $\fb\fp_2$, such that $(\alpha)=\fb\fd\fp_1\fp_2$ and $\n\fp_1\ll (\n \fb)^\ve$. 
Let us write $\fq=\fp_1\fp_2$ in what follows.
In particular we have 
$\hat{\fb}=\alpha^{-1}\fq$. Let $B_1^{(l)}=\Bl|\alpha^{(l)}|$ and $B_1=\bigoplus_l
\Blp$. In order to apply our estimates for exponential sums from the preceding section, we invoke the connection in 
\eqref{eq:call} to get
\begin{align*}
 E(v,P)
 &=
\cH(v)^{-n/2}
\sum_{
Q_0\ll \n \mathfrak{b}  \ll Q_1} 
(\n \mathfrak{b}  )^{-n}
\sum_{
\substack{\0\neq \m\in (\alpha^{-1}\fq)^n\\
|\m^{(l)}|\ll B^{(l)}}}
|\tilde{S}_\mathfrak{b}  (\alpha \m)|\\
&=\cH(v)^{-n/2}\sum\limits_{Q_0\ll\n\bbb \ll Q_1} (\n \bbb )^{-n}
\sum_{\substack{\ma{0}\neq \m\in \fq^n\\
| \m^{(l)}| \ll \Blp}} |\tilde S_\bbb(\m)|.
\end{align*}
Observe that  $\nm(B)\ll P^{-d+\ve}\cH(v)$ and $\nm(\alpha)\ll \n{\fb}\n \fq$, whence
\begin{equation}\label{lem:nm-B1}
\nm(B_1)
 \ll (\n \fb)^{1+\ve}P^{-d+\ve}\cH(v).
\end{equation}
In all that follows we will  use the  notation $q=\n\bbb$, $q_1=\n\bbb_1 $ and
$q_2=\n\bbb_2$.

Lemma \ref{lem:hb-13} implies that $\tilde S_{\bbb_1}(\m)\ll q_1^{(n+1)/2+\ve}(\n \fh)^{1/2}$, where $\fh$ is the greatest common divisor of $\fb_1$ and $(G(\m))$.
Hence, since $\fq\subseteq \fp_2$,  it follows from the multiplicativity property \eqref{eq:Lucy}  that 
\begin{align*}
 E(v,P)
&\ll P^\ve \cH(v)^{-n/2}\sum\limits_{\substack{Q_0\ll\n\bbb \ll Q_1\\ \fb=\fb_1\fb_2}} q^{-n}q_1^{(n+1)/2}
\sum_{\substack{\ma{0}\neq \m\in \fp_2^n\\
| \m^{(l)}| \ll \Blp}}  (\n \fh)^{1/2} |\tilde S_{\bbb_2}(\m)|.
\end{align*}
Our goal is to show that \eqref{eq:white} holds for a suitable absolute constant $\Delta>0$. 
On breaking the sum over $\fb_1,\fb_2$ into dyadic intervals for $\n \fb_1$ and $\n \fb_2$, it will suffice to show that 
$$
\int_R 
\max_{\substack{M_1,M_2\gg 1\\
Q_0\ll M_1M_2\ll Q_1}}
E(v,P;\M) \d v \ll P^{-\Delta},
$$
where $E(v,P;\M)$ denotes the contribution to the right hand side of the above estimate for $E(v,P)$ from $\fb_1, \fb_2$ such that $M_i\leq q_i< 2M_i$, for $i=1,2$.

It is now time to distinguish between whether $G(\m)$ is zero or non-zero in the summation over $\m$, where $G$ is the dual form that we met in \S \ref{s:sums}.
Accordingly, we  write 
$$
E(v,P;\M)= 
P^\ve \cH(v)^{-n/2}\left(
E_1(v,P;\M)+E_2(v,P;\M)\right),
$$
where
\begin{align*}
E_1(v,P;\M)
&=\sum_{\substack{\fb_1,\fb_2\\
M_i\leq q_i<2M_i}}
q^{-n}q_1^{(n+1)/2}
\sum_{\substack{\m\in \fp_2^n\\
G(\m)\neq 0\\
| \m^{(l)}| \ll \Blp}}  (\n \fh)^{1/2} |\tilde S_{\bbb_2}(\m)|,\\
E_2(v,P;\M)
&= \sum_{\substack{\fb_1,\fb_2\\
M_i\leq q_i<2M_i}}
q^{-n}q_1^{n/2+1}
\sum_{\substack{\ma{0}\neq \m\in \fp_2^n\\
G(\m)=0\\
| \m^{(l)}| \ll \Blp}} |\tilde S_{\bbb_2}(\m)|.
\end{align*}
It now suffices to show that 
\begin{equation}\label{eq:white'}
F_i=\int_R  \cH(v)^{-n/2}
\max_{\substack{M_1,M_2\gg 1\\
Q_0\ll M_1M_2\ll Q_1}}
E_i(v,P;\mathbf{M}) \d v \ll P^{-\Delta},
\end{equation}
for $i=1,2$, for a suitable absolute constant $\Delta>0$. 
Let us put  $M=M_1M_2$.

\subsection{Treatment of $F_1$}

Let $A\in \fo$ be non-zero
and let $(\fb_1,A)$ denote the greatest common divisor of $\fb_1$ and $(A)$.
Then there are at most $O(A^\ve)$ ideal divisors of $A$ and it follows that 
\begin{align*}
\sum_{\substack{\fb_1\\ 
M_1\leq q_1< 2M_1}}
\left(\n(\bbb_1,A)\right)^{1/2}
&\leq
\sum_{\fc\mid (A)}(\n\fc)^{1/2}\sum_{\substack{\fc\mid \bbb_1\\ M_1\leq q_1< 2M_1}} 
1\\
&\ll \sum_{\fc\mid (A)}(\n\fc)^{1/2}\left(\frac{M_1}{\n\fc}\right)\\
&\ll M_1 A^{\ve}. 
\end{align*}
Applying this with $A=G(\m)$, we deduce that 
\begin{align*}
E_1(v,P;\M)
&\ll P^\ve
M_1^{(3-n)/2}M_2^{-n}
\sum_{\substack{\fb_2\\
M_2\leq q_2<2M_2}}
\sum_{\substack{\ma{0}\neq \m\in \fp_2^n\\
| \m^{(l)}| \ll \Blp}}  |\tilde S_{\bbb_2}(\m)|.
\end{align*}
Now it is clear that there are $O(M_2^{1/2+\ve})$ square-full integral ideals $\fb_2$ with norm of order $M_2$. 
 Lemma \ref{lem:skinner-21}  and \eqref{lem:nm-B1} therefore yield
\begin{align*}
E_1(v,P;\M)
&\ll P^\ve
M_1^{(3-n)/2}M_2^{-n}
\sum_{\substack{\fb_2\\
M_2\leq q_2<2M_2}}
M_2^{n/2+1}(\nm(B_1)^n+M_2^{n/3})\\
&\ll P^\ve
M_1^{(3-n)/2}M_2^{-n}
M_2^{(n+3)/2}( P^{-dn}\cH(v)^n M^n+M_2^{n/3})\\
&\ll P^\ve
M^{(3-n)/2}
( P^{-dn}\cH(v)^n M^n+M^{n/3})\\
&= P^\ve
( P^{-dn}\cH(v)^n M^{(n+3)/2}+M^{3/2-n/6}).
\end{align*}
Notice that the exponent $3/2-n/6$ is negative for $n\geq 10$. 
Since $Q_0\ll M\ll Q_1$ and $n\geq 10$, we
therefore obtain
\begin{align*}
E_1(v,P;\M)
&\ll P^\ve
\left( P^{-dn}\cH(v)^n Q_1^{(n+3)/2}+Q_0^{3/2-n/6}\right)\\
&= P^\ve
\left( P^{-d(n-9)/4}\cH(v)^{(n-3)/2}+P^{-d(n-9)/6}\cH(v)^{n/6-3/2}\right)\\
&\ll P^{-d/6+\ve}
\left( \cH(v)^{(n-3)/2}+\cH(v)^{n/6-3/2}\right).
\end{align*}
Inserting this into \eqref{eq:white'} and applying Lemma \ref{nmintegral}, this therefore shows that 
\begin{align*}
F_1\
&\ll P^{-d/6+\ve}
\int_R \cH(v)^{-3/2}\d v
\ll P^{-d/6+\ve},
\end{align*}
which is satisfactory.

\subsection{Treatment of $F_2$}

According to  Lemmas \ref{lem:skinner-21} and \ref{lem:hb-16} we 
have 
\begin{align*}
E_2(v,P;\M)
&\ll  \sum_{\substack{\fb_1,\fb_2\\
M_i\leq q_i<2M_i}}
M^{-n}M_1^{n/2+1}
\sum_{\substack{\ma{0}\neq \m\in \fp_2^n\\
G(\m)=0\\
| \m^{(l)}| \ll \Blp}} |\tilde S_{\bbb_2}(\m)|\\
&\ll P^\ve  \sum_{\substack{\fb_1,\fb_2\\
M_i\leq q_i<2M_i}}
M^{-n}M_1^{n/2+1}
M_2^{n/2+1} \mathcal{M}.
\end{align*}
Here
\begin{align*}
\mathcal{M}
&=\min\left\{
 \nm(B_1)^n +M_2^{n/3}, 
M_2^{(n-1)/2}+\nm(B_1)^{n-3/2} M_2^{1/3}\right\}\\
&\leq M_2^{n/3}+ \nm(B_1)^{n-3/2} M_2^{1/3}+
 \min\left\{
\nm(B_1)^n,
M_2^{(n-1)/2}\right\}\\
&=\mathcal{M}_1+\mathcal{M}_2+\mathcal{M}_3,
\end{align*}
say.
We write 
$E_{2,j}$ for the overall contribution to 
$E_2(v,P;\mathbf{M})
$ from  $\mathcal{M}_j$, for $j\in \{1,2,3\}$.

Beginning with the contribution from $\cM_1$, we sum over $O(M_1)$ ideals 
$\fb_1$ and $O(M_2^{1/2+\ve})$ square-full ideals $\fb_2$, finding that 
\begin{align*}
E_{2,1}
&\ll P^\ve  \sum_{\substack{\fb_1,\fb_2\\
M_i\leq q_i<2M_i}}
M^{-n}M_1^{n/2+1}
M_2^{5n/6+1} \\
&\ll P^\ve  
M_1^{2-n/2}
M_2^{3/2-n/6}\\ 
&\ll P^\ve  
M^{(9-n)/6},
\end{align*}
since $2-n/2<3/2-n/6<0$  for $n\geq 10$.
Taking $M\gg Q_0$, inserting the outcome  into 
 the left hand side of 
\eqref{eq:white'} and then applying Lemma \ref{nmintegral}, this therefore shows that this case  contributes
$O( P^{-d/6+\ve})$
to $F_2$, 
which is satisfactory.

Turning to  the contribution from  $\cM_2$, we apply \eqref{lem:nm-B1} and 
sum over $\fb_1$ and $\fb_2$, as before. In this way  we obtain
\begin{align*}
E_{2,2}
&\ll P^\ve  \sum_{\substack{\fb_1,\fb_2\\
M_i\leq q_i<2M_i}}
M^{-n}M_1^{n/2+1}
M_2^{n/2+1} 
\nm(B_1)^{n-3/2} M_2^{1/3}\\
&\ll P^{-d(n-3/2)+\ve} \cH(v)^{n-3/2}  \sum_{\substack{\fb_1,\fb_2\\
M_i\leq q_i<2M_i}}
M^{-3/2}M_1^{n/2+1}
M_2^{n/2+4/3} \\
&\ll P^{-d(n-3/2)+\ve} \cH(v)^{n-3/2}  
M_1^{n/2+1/2}
M_2^{n/2+1/3} \\
&\ll P^{-d(n-3/2)+\ve} \cH(v)^{n-3/2}  
M^{n/2+1/2}.
\end{align*}
Taking $M\ll Q_1$, inserting the outcome  into
 the left hand side of  \eqref{eq:white'} and then applying Lemma \ref{nmintegral}, the $\cM_2$ term is therefore seen to contribute
\begin{align*}
&\ll
P^{-d(n-3/2)+3d(n+1)/4+\ve} \int_R \cH(v)^{-2}   \d v\\
&\ll   
P^{d(-n/4 +9/4)+\ve} \\
&\ll   
P^{-d/4+\ve}
\end{align*}
to $F_2$, 
since $n\geq 10$.

It remains to deal with the term $\cM_3$, for which we take 
$\min\{A,B\}\leq A^{(n-2)/(n-1)}B^{1/(n-1)}$. In view of \eqref{lem:nm-B1}, this gives
$$
\cM_3\ll P^\ve
\left(
M P^{-d} \cH(v)\right)^{n(n-2)/(n-1)}
M_2^{1/2}.
$$
Summing over $\fb_1,\fb_2$ as before, we 
obtain 
\begin{align*}
E_{2,3}
&\ll P^\ve \left(
P^{-d} \cH(v)\right)^{n(n-2)/(n-1)}
  \sum_{\substack{\fb_1,\fb_2\\
M_i\leq q_i<2M_i}}
M^{-n/(n-1)}M_1^{n/2+1}
M_2^{n/2+3/2}\\
&\ll P^\ve \left(
P^{-d} \cH(v)\right)^{n(n-2)/(n-1)}
M^{n/2+2-n/(n-1)}.
\end{align*}
Taking $M\ll Q_1$, inserting the outcome  into the left hand side of  \eqref{eq:white'} and applying Lemma~\ref{nmintegral}, the $\cM_3$ term  therefore contributes
\begin{align*}
&\ll
P^{-dn(n-2)/(n-1)+\ve}
\int_R \cH(v)^{-n/2+
n(n-2)/(n-1)}
Q_1^{n/2+2-n/(n-1)} \d v\\
&\ll
P^{-dn(n-2)/(n-1)+3d(n/4+1-n/(2n-2))+\ve}
\int_R \cH(v)^{-2} \d v\\
&\ll
P^{-dn(n-2)/(n-1)+3d(n/4+1-n/(2n-2))+\ve}
\end{align*}
to $F_2$.  This is 
$O(P^{-d/18+\ve})$, since 
 $n\geq 10$, which  therefore concludes the proof that \eqref{eq:white'} holds for $i=1,2$.

\subsection{Conclusion}

Bringing everything together, our work so far has established the following result.

\begin{lemma}
\label{main sum}
Assume that $n\geq 10$.
Then there exists an absolute constant $\Delta>0$ such that 
$$
\No=
\frac{c_{Q}2^{r_2 n}}{D_{K}^{n/2}Q^{2d}}
\sum_{\substack{(0)\neq \mathfrak{b}\subseteq \fo\\
\n \mathfrak{b}  \ll Q^{d}
}} (\n \mathfrak{b}  )^{-n}S_\mathfrak{b}  (\0)
I_\mathfrak{b}  (\0) 
+
O(P^{d(n-3)-\Delta}).
$$
\end{lemma}

Our remaining task is to show that the main term here converges to the main term predicted in Theorem \ref{t:10'}.
Recall from \S \ref{s:integrals} that $\rho=Q^{-d}\n\bbb$. We therefore proceed to 
consider the sum
$$
\sum_{\substack{(0)\neq \mathfrak{b}\subseteq \fo\\
\n \mathfrak{b}  \ll Q^{d}
}} (\n \mathfrak{b}  )^{-n}S_\mathfrak{b}  (\0)
I_\fb  (\0).
$$ 
According to \eqref{eq:def-I} and \eqref{eq:exp-int'} we have 
\begin{align*}
I_\fb(\0)
&=P^{dn}\int_{V^n} W(\x)h(\rho,\nm(F(\x)))\d \x\\
&=P^{dn}\int_V \wtil(v)h(\rho,\nm(v))\d v.
\end{align*}
Now it follows from Lemma \ref{prop:1} that 
$$
\int_V
f(v)h(\rho,\nm(v))\d v=\frac{\sqrt{D_K}}{2^{r_2}} f(0)+O_N(\lambda_f^{2d(N+1)}\rho^N),
$$
for any $N>0$ and any $f\in \cW_1(V)$, in the notation of \eqref{eq:sobolev}.
Recalling the properties of the weight function $\wtil$ discussed in \S \ref{s:poisson-forms},
and in particular \eqref{eq:r1}, it therefore follows that 
\begin{equation}
I_\fb(\0)
=\frac{\sqrt{D_K}}{2^{r_2}} P^{dn}\wtil(0)+ O_N\left(
P^{dn+\ve} \rho^{N} 
\right).
\label{eq:singular integral}
\end{equation}
The analysis of $\wtil(0)$ is the object of the following result.

\begin{lemma}\label{lem:I0}
We have $(\log P)^{2d(1-n)}\ll\wtil(0)\ll (\log P)^{2d(1-n)}.$
\end{lemma}

\begin{proof}
Recall that $W(\x)=w_0(\x)\omega(\x),$ with $w_0$ and $\omega$ as in \S \ref{s:poisson-forms}. 
In particular,  outside the set $\|\x^{(l)}-\bxi^{(l)}\| \leq
(\log P)^{-1}$ we have  $\omega(\x^{(l)}) =O(e^{-(\log P)^2})$, whereas  $w_0(\x^{(l)})=1$ within it, provided that $P$ is sufficiently large.
It follows  from \eqref{eq:black} that 
\begin{align*}
 \wtil(0)
&=\prod_{l=1}^{r_1+r_2}
I^{(l)}(0)
+O(e^{-(\log
P)^2}),
\end{align*}
where
$$
I^{(l)}(0)=
\int \omega(\x^{(l)})\d\x^{(l)},
$$
the integral 
being taken over $\x^{(l)}\in K_l^n$ for which 
$F^{(l)}(\x^{(l)})=0$ and $\|
\x^{(l)}-\bxi^{(l)}\|\leq (\log P)^{-1}$.

We begin by analysing $I^{(l)}(0)$ for $l\leq r_1$, making the change of variables $\u=\x-\bxi$ and then 
$$
v_1=F^{(l)}(\u^{(l)}+\bxi^{(l)}),
\quad \mbox{$v_i=u^{(l)}_i$ for $2\leq i\leq n$}.
$$
According to \eqref{eq:space} we have $u^{(l)}_1=f^{(l)}(\v)$, whence $\d\x^{(l)}=\d\u^{(l)}=
\frac{\partial f^{(l)}}{\partial v_1}\d\v$. 
Let 
$$
g(v)=f^{(l)}(\v)^2+v_2^2+\cdots +v_n^2.
$$
It now follows that 
\begin{align*}
I^{(l)}(0)&=
J^{(l)}(P)
+O(e^{-(\log P)^2}),
\end{align*}
with 
\begin{align*}
J^{(l)}(P)=
\int \exp\left(-(\log
P)^4 g(0,v_2,\ldots,v_n)\right)\frac{\partial f^{(l)}(0,v_2,...,v_n)}{\partial
v_1}\d v_2\cdots \d v_n,
\end{align*}
where the integral is over 
$v_2,\ldots,v_n\in \RR$ such that 
$g(0,v_2,\ldots,v_n)\leq (\log P)^{-2}$. 
The term 
 $J^{(l)}(P)$ precisely coincides with the term $I_1(0)$ from \cite[page 252]{hb-10},  the analysis of which shows that 
$(\log P)^{2(1-n)}\ll I^{(l)}(0)\ll (\log P)^{2(1-n)} $. 

An analogous method can be used to show that  when $l>r_1$, one has
$$
(\log P)^{4(1-n)}\ll I^{(l)}(0)\ll (\log P)^{4(1-n)}.
$$ 
In fact one finds that $I^{(l)}(0)$ is asymptotically equal to $I_2(0,0)$, in the notation of  
\cite[\S 12]{skinner}.
Taken together with our work above, this leads to the conclusion of the lemma.
\end{proof}

Taking $N=1$ in  \eqref{eq:singular integral}
we obtain
\begin{equation}\label{eq:skool}
\begin{split}
\sum_{\substack{
\n \mathfrak{b}  \ll Q^{d}
}} (\n \mathfrak{b}  )^{-n}S_\mathfrak{b}  (\0)
I_\mathfrak{b}  (\0) 
=~&
\frac{\sqrt{D_K}}{2^{r_2}}
P^{dn}\mathfrak{I}(0)
\sum_{
\n \mathfrak{b}  \ll Q^{d}
} (\n \mathfrak{b}  )^{-n}S_\mathfrak{b}  (\0)\\
&+
O\left(P^{dn+\ve}
\sum_{
\n \mathfrak{b}  \ll Q^{d}
}  \rho(\n \mathfrak{b}  )^{-n}|S_\mathfrak{b}  (\0)|
\right).
\end{split}
\end{equation}
We begin by dealing with the  error term in this expression. 
Thus let  
$$
\fS_1=\sum_\bbb
(\n\bbb)^{-n}|S_\bbb(\0)|,
$$ 
the sum being over all non-zero integral ideals $\fb$.
Recalling that $\rho=Q^{-d}\n\fb$,  we see that
\begin{align*}
\sum_{\n\bbb\ll Q^d}\rho (\n \mathfrak{b}  )^{-n}|S_\mathfrak{b}(\0)|&\ll P^{-1}\sum_{\n\bbb\ll
P^{-1}Q^d}(\n
\mathfrak{b}  )^{-n}|S_\mathfrak{b}(\0)|+\sum_{\n\bbb\gg P^{-1}Q^d}(\n \mathfrak{b} 
)^{-n}|S_\mathfrak{b}(\0)|\\
&\ll P^{-1}\fS_1+\sum_{\n\bbb\gg P^{-1}Q^d}(\n \mathfrak{b}  )^{-n}|S_\mathfrak{b}(\0)|.
\end{align*}
Adopting the notation $q,q_1,q_2$ from above, it follows from Lemmas \ref{lem:hb-13} and \ref{lem:skinner-21} that
$$
S_\bbb(\0)\ll q_1^{n/2+1+\ve }q_2^{5n/6+1+\ve}.
$$ 
Hence, on breaking into dyadic intervals for $q_1,q_2$,  we easily deduce that  
\begin{align*}
\sum_{\n\bbb\geq A} (\n \mathfrak{b}  )^{-n}|S_\mathfrak{b}(\0)|
&\ll
\sum_{\substack{\fb_1,\fb_2\\ 
q\geq A}}
q_1^{1-n/2+\ve}q_2^{1-n/6+\ve}\ll A^{3/2-n/6+\ve}\ll A^{-1/6+\ve},
\end{align*}
for  $n\geq 10$.
In particular this establishes the convergence of  $\fS_1$ and, a fortiori, the absolute convergence of 
$\fS$ in \eqref{eq:sing-series}.
We deduce that 
$$
\sum_{\substack{
\n \mathfrak{b}  \ll Q^{d}
}} (\n \mathfrak{b}  )^{-n}S_\mathfrak{b}  (\0)= \fS+ O(P^{-3/12+\ve})
$$
and 
\begin{align*}
\sum_{\n\bbb\ll Q^d}\rho (\n \mathfrak{b}  )^{-n}|S_\mathfrak{b}(\0)|
&\ll P^{-1} +(P^{-1}Q^d)^{-1/6+\ve}
\ll P^{-1/12+\ve}, 
\end{align*}
in \eqref{eq:skool}.
Once inserted into Lemma
\ref{main sum}, under the assumption that $n\geq 10$, we see that there is an absolute constant $\Delta>0$ such that 
$$
\No=
\frac{c_{Q}2^{r_2 (n-1)}}{D_{K}^{(n-1)/2}} \fS \wtil(0)
P^{d(n-3)} 
+
O(P^{d(n-3)-\Delta}).
$$
Here $c_Q=1+O_N(Q^{-N})$. Moreover, 
the  lower bound $\fS\gg 1$ is standard and can be established using the argument of Pleasants  \cite[Lemma~7.4]{pleasants}, for example. 
Hence an application of Lemma \ref{lem:I0} shows that 
$\fS \wtil(0)$ has the order of magnitude claimed in 
Theorem \ref{t:10'}, thereby concluding our argument.

\end{document}